\documentclass[12pt,reqno]{amsart}
\input{epsf}
\usepackage{amsthm,amsmath,amssymb,amsfonts,amscd,verbatim,latexsym,xypic,graphics}

\newtheorem{Theorem}{Theorem}[section]
\newtheorem{Lemma}[Theorem]{Lemma}
\newtheorem{Proposition}[Theorem]{Proposition}
\newtheorem{Corollary}[Theorem]{Corollary}

\newtheorem{Example}[Theorem]{Example}

\newcommand{\constant}{\Bbbk } 
\renewcommand{\O}{\mathcal O}
\newcommand{\M}{\mathcal M}
\newcommand{\complex}{\mathbf C}
\newcommand{\rationals}{\mathbf Q}
\newcommand{\rank}{\text{rank} \,}
\newcommand{\ux}{\mathbf x}
\newcommand{\uy}{\mathbf y}
\newcommand{\uz}{\mathbf z}
\newcommand{\hf}{\hat f}
\renewcommand{\P}{\mathbb P}
\newcommand{\cS}{\mathcal S}
\newcommand{\bU}{\mathbb U} 
\newcommand{\A}{\mathcal A}
\newcommand{\C}{\mathcal C}
\newcommand{\W}{\mathcal W}

\newcommand{\G}{\mathcal G}
\newcommand{\HI}{\mathbb H}
\newcommand{\II}{\mathbb I}
\newcommand{\I}{\mathfrak I}
\newcommand{\X}{\mathcal X}
\newcommand{\He}{\mathcal H}
\newcommand{\Sym}{\text{Sym}}
\newcommand{\Proj}{\text{Proj} \, }
\newcommand{\cb}{\,) \hspace{-1.6mm} ( \,}
\newcommand{\F}{\mathbb F}
\newcommand{\E}{\mathcal E}
\newcommand{\cF}{\mathcal F}
\newcommand{\cU}{\mathcal U}
\newcommand{\cV}{\mathcal V}
\newcommand{\gb}{\mathfrak b}
\newcommand{\R}{\mathfrak R}
\newcommand{\B}{\mathcal B}
\newcommand{\T}{\mathbb T}
\newcommand{\ssz}{\mathsf z}
\newcommand{\ssa}{\mathsf a}
\newcommand{\ssb}{\mathsf b}
\newcommand{\ssc}{\mathsf c}
\newcommand{\ssd}{\mathsf d}
\newcommand{\ssp}{\mathsf p}
\newcommand{\ssq}{\mathsf q}
\newcommand{\ssw}{\mathsf w}
\newcommand{\cZ}{\mathcal Z}

\newcommand{\fp}{\mathfrak p}
\newcommand{\fq}{\mathfrak q}
\newcommand{\en}[1]{\langle #1 \rangle}
\newcommand{\ra}{\rightarrow}

\newcommand{\lra}{\longrightarrow}
\newcommand{\la}{\leftarrow}

\newcommand{\demo}{\noindent {\sc Proof.}\;}

\begin{document} 
\title{On Hermite's invariant for binary quintics}
\author{Jaydeep Chipalkatti} 
\maketitle 

\parbox{12.5cm}{\small 
{\sc Abstract.} 
Let $\He \subseteq \P^5$ denote the hypersurface of binary quintics in involution, with 
defining equation given by the Hermite invariant $\HI$. In \S\ref{section.sing.H} we find the singular 
locus of $\He$, and show that it is a complete intersection of a 
linear covariant of quintics. In \S\ref{section.dual.H} we show that the 
projective dual of $\He$ can be canonically identified with itself via an involution. 
The Jacobian ideal of $\HI$ is shown to be perfect of height two in \S\ref{section.J}, moreover 
we describe its $SL_2$-equivariant minimal free resolution. The last section 
develops a general formalism for evectants of covariants of binary forms, which is then used to 
calculate the evectant of $\HI$.} 

\bigskip \medskip 

\parbox{12cm}{\small 
} 

\bigskip  

\section{Introduction} \label{section.introduction} 

This paper analyses the geometry and invariant theory of the 
Hermite invariant for binary quintics. We begin by recalling the 
elementary properties of this invariant; the main results are summarised 
on pages~\pageref{section.results.summary}-\pageref{end.summary} 
after the required notation is available. 
We refer to~\cite{Glenn,GrYo} and~\cite{Salmon1} for foundational notions in 
the classical invariant theory of binary forms, as well as the symbolic 
method. Modern treatments of 
this material may be found in~\cite{Dolgachev,Gurevich,Kung-Rota} 
and~\cite{Olver}. 
The encyclop{\ae}dia article~\cite{MacMahon} contains a very readable 
introduction to the classical theory. 
We will use~\cite[Lecture 11]{FH} and \cite[\S4.2]{Sturmfels} 
for the basic representation theory of $SL_2$. 
The discovery of the Hermite invariant was first reported in 
\cite[Premi{\`e}re Partie, \S IV-VII]{Hermite1}. 

\begin{figure} 
\epsffile[-70 10 187 300]{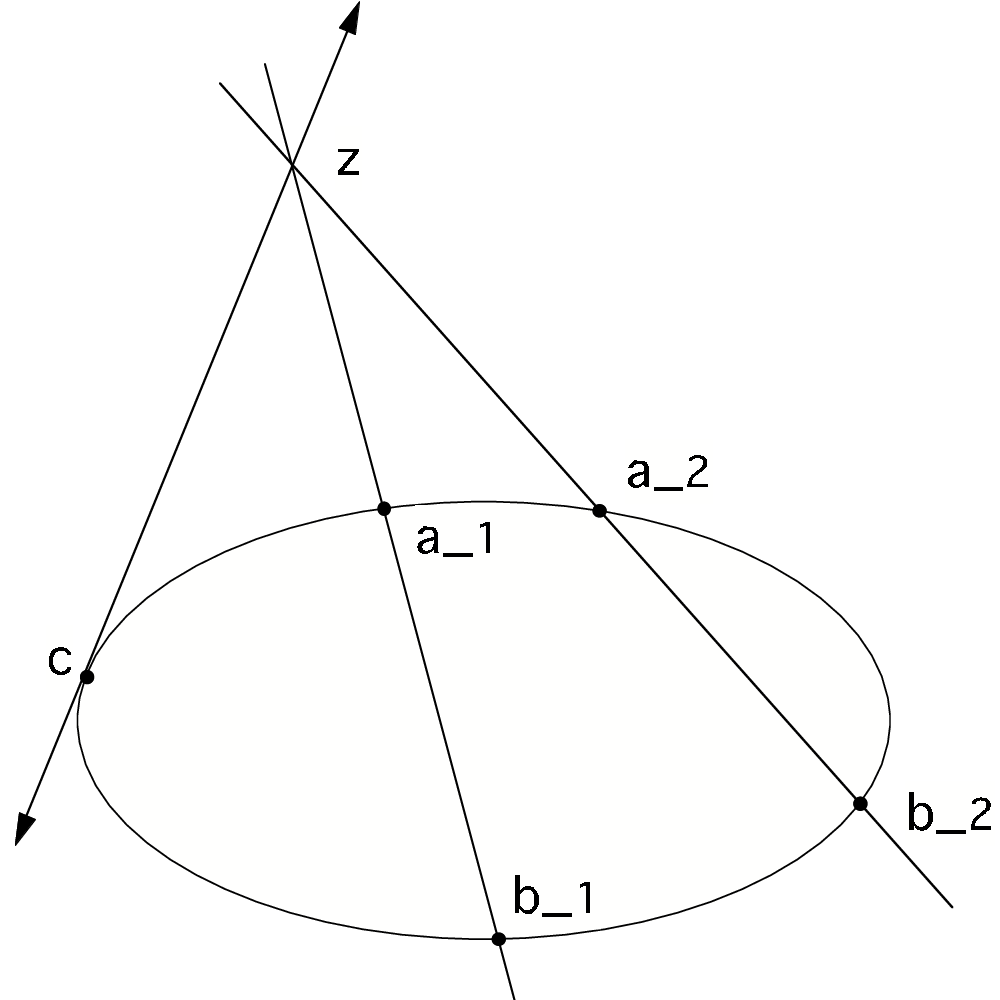} 
\end{figure} 

The results in Lemma~\ref{lemma.Hdegree} and 
Proposition~\ref{prop.triple.intersection} below are classical; I 
have included them for completeness of treatment. 
\subsection{} \label{section.conic_diagram}
The base field will be $\complex$. Let $V$ denote a two-dimensional 
complex vector space with basis $\ux = \{x_1,x_2\}$ and a natural action of 
$SL(V)$. For $m \ge 0$, let $S_m = \Sym^m \, V$ denote the 
$(m+1)$-dimensional irreducible $SL(V)$-representation consisting of 
binary $m$-ics in $\ux$. Consider the quadratic Veronese imbedding 
\[ \phi: \P \, V \lra \P S_2, \quad 
[c_1 \, x_1 + c_2 \, x_2] \lra 
[(c_1 \, x_1+c_2 \, x_2)^2], \] 
whose image is a smooth conic $\phi(\P^1) = C \subseteq \P^2$. We identify 
$\P^5$ with $\Sym^5 \, C \simeq \P S_5$, 
i.e., a point in $\P^5$ is alternately seen as 
a degree $5$ effective divisor on $C$, or as a 
binary quintic in $\ux$ distinguished up to scalars. 

\label{conic_diagram} 
Let $\ssz$ be a point of $\P^2 \setminus C$, and $L_1,L_2$ two lines through 
$\ssz$ intersecting $C$ in $\ssa_1,\ssb_1;\ssa_2,\ssb_2$. Let $\ssc \in C$ be 
one of the two points such that the line $\overline{\ssc \, \ssz}$ 
is tangent to $C$, and now define a divisor
$\ssa_1 + \ssb_1 + \ssa_2 + \ssb_2 + \ssc \in \P^5$. As $\ssz,L_1,L_2$ move, let 
$\He \subseteq \P^5$ denote the closure of the set of all such divisors. 
(The closure includes all divisors of the form $3 \, \ssz + \ssa + \ssb$ for arbitrary 
points $\ssz,\ssa,\ssb$ in $C$.) 

There are $\infty^2$ possible positions for $\ssz$, and then 
$\infty^1$ positions for each of the $L_i$ once $\ssz$ is fixed; 
hence $\dim \He = 4$. By construction $\He$ is an irreducible variety. 
The action of $SL(V)$ on $\P S_2$ induces an action 
on $C$, moreover it takes a tangent line
to $C$ to another tangent line, hence $SL(V)$ acts on the imbedding
$\He \subseteq \P^5$.  Consequently the equation of $\He$ is an invariant 
of binary quintics, usually called the Hermite invariant $\HI$. This defines 
$\HI$ only up to a multiplicative constant; but see formula~(\ref{defn.H}) below. 

A point $\ssz \in \P^2 \setminus C$ defines an order 
$2$ automorphism of $C$, sending $\ssa \in C$ to the other intersection 
of $\overline{\ssz \, \ssa}$ with $C$. 
The divisor $\ssz + \ssa_1 + \ssb_1 + \ssa_2 + \ssb_2$ is said to be 
{\sl in involution} with respect to $\ssz$ 
since it is fixed by this automorphism. 

\begin{Lemma} \sl The degree of $\He$ is $18$. 
\label{lemma.Hdegree} \end{Lemma} 
\demo 
For $\ssp \in C$, let $\Gamma_\ssp \subseteq \P^5$ denote the hyperplane 
defined by all the divisors containing $\ssp$. Given general points 
$\ssp_1,\ssp_2,\ssp_3,\ssp_4$ in $C$, consider the intersection 
$\Sigma = \He \cap \Gamma_{\ssp_1} \cap \dots \cap \Gamma_{\ssp_4}$. The 
three points 
\begin{equation} \overline{\ssp_1 \, \ssp_2} \cap \overline{\ssp_3 \, \ssp_4}, 
\quad 
\overline{\ssp_1 \, \ssp_3} \cap \overline{\ssp_2 \, \ssp_4}, \quad 
\overline{\ssp_1 \, \ssp_4} \cap \overline{\ssp_2 \, \ssp_3}, 
\label{triple-intersection} \end{equation}
give $6$ elements in $\Sigma$ (since two tangents to $C$ can be drawn 
from each). Alternately, let the tangent to $C$ at $\ssp_1$ intersect 
$\overline{\ssp_2 \, \ssp_3}$ at $\ssz$, and let 
$\overline{\ssz \, \ssp_4}$ intersect $C$ in the additional point $\ssq$; which 
gives $\ssp_1 + \dots + \ssp_4 + \ssq \in \Sigma$. This construction produces 
$4 \times 3 = 12$ more elements in $\Sigma$, hence 
$\text{card} \, (\Sigma) = 18$. \qed 

\smallskip 

\subsection{} \label{section.FQ}
With notation as in the diagram, write $\ssc = [\phi(x_1)]$ after a change of variables. Then 
$\ssa_1,\ssb_1$ must equal $\phi([\alpha_1 x_1 + \alpha_2 x_2]),
\phi([\alpha_1 x_1 - \alpha_2 x_2])$ for some $[\alpha_1,\alpha_2] \in \P^1$, and 
similarly for $\ssa_2,\ssb_2$. Hence $\ssa_1 + \ssa_2 + \ssb_1 + \ssb_2 + \ssc$ corresponds to the quintic 
\begin{equation} \cF_Q = x_1 \, (q_0 \, x_1^4 + 2 \, q_1 \, x_1^2 \, x_2^2 + q_2 \, x_2^4)
\end{equation} 
for some $Q=[q_0,q_1,q_2] \in \P^2$. This `canonical form' will prove most useful for computations. 
Since any $[F] \in \He$ lies in the $SL_2$-orbit of some $[\cF_Q]$, any `equivariant' calculation which is valid for 
$\cF_Q$ is valid generally. 

\smallskip 

In the next few sections we will gather some needed preliminaries from 
classical invariant theory; we will take up $\HI$ once more on page~\pageref{defn.H}. 

\subsection{Transvectants} \label{section.trans} 
Given integers $m,n \ge 0$, we have a decomposition of 
$SL(V)$-representations 
\begin{equation} 
S_m \otimes S_n \simeq \bigoplus\limits_{r=0}^{\min(m,n)} \, 
S_{m+n-2r}.   \label{Clebsch-Gordan} \end{equation}
Let $A,B$ denote binary forms in $\ux$ of respective orders $m,n$. The $r$-th transvectant 
of $A$ with $B$, written $(A,B)_r$, is defined to be the image of 
$A \otimes B$ via the projection map 
\[ \pi_r: S_m \otimes S_n \lra S_{m+n-2r} \, .  \] 
It is given by the formula 
\begin{equation} (A,B)_r = \frac{(m-r)! \, (n-r)!}{m! \, n!} \, 
\sum\limits_{i=0}^r \, (-1)^i \binom{r}{i} \, 
\frac{\partial^r A}{\partial x_1^{r-i} \, \partial x_2^i} \, 
\frac{\partial^r B}{\partial x_1^i \, \partial x_2^{r-i}} 
\label{trans.formula} \end{equation} 
(Some authors choose the initial 
scaling factor differently, cf.~\cite[Ch.~5]{Olver}.) 
By convention $(A,B)_r = 0$ if $r > \min \, (m,n)$. If we symbolically write 
$A = \alpha_\ux^m, B = \beta_\ux^n$, then 
$(A,B)_r = (\alpha \, \beta)^r \, \alpha_\ux^{m-r} \, \beta_\ux^{n-r}$. 
There is a canonical isomorphism of representations 
\begin{equation} 
 S_m \stackrel{\sim}{\lra} S_m^* \, ( \, = \text{Hom}_{SL(V)}(S_m,S_0)) 
\label{self-duality} \end{equation} 
which sends $A \in S_m$ to the functional $B \lra (A,B)_m$.  Hence 
if $A$ is an order $m$ form such that $(A,B)_m=0$ for all $B \in S_m$, then 
$A$ must be zero. 
\subsection{Gordan series} \label{section.Gordanseries} 
Introduce a parallel set of letters 
$\uy = (y_1,y_2)$, and define Cayley's Omega operator 
\[ \Omega_{\ux \uy} = 
\frac{\partial^2}{\partial x_1 \, \partial y_2} - 
\frac{\partial^2}{\partial x_2 \, \partial y_1}. \] 
If we represent an element in $S_m \otimes S_n$ as a bihomogeneous form $G$
of orders $m,n$ in $\ux,\uy$, then \[ \pi_r (G) = \frac{(m-r)! \, (n-r)!}{m! \, n!} \, 
\{\Omega_{\ux \uy}^r \circ G \}_{\uy:=\ux \, .} \] 
A splitting to $\pi_r$ is given by the map 
\[ \imath_r: \alpha_{\ux}^{m+n-2r} \lra (\ux \, \uy)^r \, 
\alpha_{\ux}^{m-r} \alpha_{\uy}^{n-r}, \] where 
$(\ux \uy) = x_1 \, y_2 - x_2 \, y_1$. The decomposition $G = \sum\limits_r \, \imath_r \circ \pi_r(G)$ is called 
the Gordan series for $G$. In general, it may be symbolically written as 
\[ \alpha_\ux^m \, \beta_\uy^n = 
\sum\limits_{r=0}^{\min(m,n)} \, 
\frac{\binom{m}{r} \, \binom{n}{r}}{\binom{m+n-r+1}{r}} \, 
(\ux \, \uy)^r \, {\theta_{(r)}}_{\ux}^{m-r} \, {\theta_{(r)}}_{\uy}^{n-r}, \] 
where ${\theta_{(r)}}_\ux^{m+n-2r}$ stands for $(\alpha \, \beta)^r \, 
\alpha_\ux^{m-r} \, \beta_\ux^{n-r}$ (see~\cite[p.~55]{GrYo} or~\cite[\S 24.4]{Gurevich}).

\subsection{Wronskians} 
Let $m,n \ge 0$ be integers such that $m \le n+1$. 
Consider the following composite morphism of representations 
\[ w: \wedge^m S_n \stackrel{\sim}{\lra} S_m(S_{n-m+1}) \lra S_{m(n-m+1)}, \] 
where the first map is an isomorphism (see~\cite[\S2.5]{AC}) and the second is the natural surjection. 

Given a sequence of binary $n$-ics $A_1,\dots,A_m$, define 
their Wronskian $W(A_1,\dots,A_m)$ to be the 
determinant 
\[ (i,j) \lra \frac{\partial^{m-1}  \, A_i}{\partial x_1^{m-j} \, \partial \, x_2^{j-1}}, \quad 
(1 \le i,j \le m). \] 
It equals the image $w(A_1 \wedge \dots \wedge A_m)$. We have 
$W(A_1,\dots,A_m)=0$, iff the $A_i$ are linearly dependent over $\complex$. (The `if' part is obvious. For 
the converse, see~\cite[\S 1.1]{Meulien}.) 

\begin{Lemma} \sl 
Let $A_1,\dots,A_m$ be linearly independent forms of order $m$. Then 
$W = W(A_1,\dots,A_m)$ is (up to scalar) the unique form of order $m$ 
such that $(W,A_i)_m=0$ for all $i$. 
\label{lemma.wr} \end{Lemma} 
\demo 
Consider the composite morphism 
\[ g: \wedge^{m+1} S_m \lra \wedge^m S_m \otimes S_m 
\stackrel{\sim}{\lra} S_m \otimes S_m \lra \complex,  \] 
where the first map is dual to the exterior product. 
For any $i$, we have $(W,A_i)_m = g(A_1 \wedge \dots \wedge A_m \wedge A_i) = 0$. 
The pairing 
\[ S_m \times S_m \lra \complex, \quad (A,B)\lra (A,B)_m \] 
is nondegenerate, hence such a form is unique up to scalar. \qed 

\subsection{Covariants} 
Reviving an old notation due to Cayley, 
we will write $(\alpha_0,\dots,\alpha_n \cb u,v)^n$ for the expression 
\[\sum\limits_{i=0}^n \; \binom{n}{i} \, \alpha_i \, u^{n-i} v^i. \] 
In particular $\F = (a_0,\dots,a_d \cb x_1,x_2)^d$ denotes the {\sl generic} $d$-ic, 
which we identify with the natural trace form in $S_d \, \otimes \, S_d^*$. 
Using the duality in~(\ref{self-duality}), this amounts to the identification 
of $a_i \in S_d^*$ with $\frac{1}{d!} \, x_2^{d-i} \, (-x_1)^i$. 
Let $R$ denote the symmetric algebra 
\[ \bigoplus\limits_{m \ge 0} \, S_m(S_d^*) = \bigoplus\limits_{m \ge 0} \, R_m = 
\complex \, [a_0,\dots,a_d], \] 
and $\P^d = \P \, S_d = \Proj \, R$. 

A {\sl covariant} of degree-order $(m,q)$ (of binary $d$-ics) is  by definition an $SL(V)$-equivariant imbedding 
$S_0 \hookrightarrow S_m(S_d) \otimes S_q$. 
Let $\Phi$ denote the image of $1$ via this map, then we may write 
$\Phi = (\varphi_0,\dots,\varphi_q \cb x_1,x_2)^q$
where each $\varphi_i$ is a homogeneous degree $m$ form in the 
$\{a_i\}$. The weight of $\Phi$ is defined to be $\frac{1}{2}(d \, m-q)$ 
(which is always a nonnegative integer). 
A covariant of order $0$ is called an invariant. 
E.g., $(\F,\F)_2$ is a covariant of degree-order $(2,2d-4)$, and for $d=4$, 
the compound transvectant $((\F,\F)_2,\F)_4$ is an invariant of degree $3$. If 
$\F$ is specialized to $F \in S_d$, then $\Phi$ gets specialized to $\Phi_F \in S_q$. 

\subsection{} \label{section.map.Rmodules}
Let $\Phi$ denote a covariant of degree-order $(m,q)$. Let $a,b$ denote 
nonnegative integers, and let $r= (a+q-b)/2$. For every $F \in S_d$,  we have a map 
\[ h_F: S_a \lra S_b, \quad G \lra (\Phi_F,G)_r. \] 
Since the entries of the matrix describing $h_F$ are degree $m$ forms in the $\{a_i\}$, we 
may see it as an $SL_2$-equivariant  map of graded $R$-modules 
\begin{equation} R \otimes S_a \lra R(m) \otimes S_b, 
\label{map.Rmodules} \end{equation} 
Conversely, every equivariant map of the form (\ref{map.Rmodules}) 
arises from a covariant. (Indeed, in degree zero it reduces to 
a map of representations $S_a \lra S_m(S_d) \otimes S_b$.) 
The numerical conditions are assumed to be such that the transvection is possible, 
i.e., we must have $a+q-b$ nonnegative and even, and 
$r \le \min(a,q)$. 

If $a \le b$, then by the Wronskian of the map $h$ we mean 
\[ W(h_\F(x_1^a),h_\F(x_1^{a-1} \, x_2), \dots, h_\F(x_2^a)), \] 
which is a covariant of degree $m \, (a+1)$ and order $(a+1)(b-a)$. Its coefficients are 
(up to signs) the maximal minors of $h_\F$. 
\subsection{} \label{section.evectant} 
We will let $\I(\Phi) \subseteq R$ denote the ideal generated by the 
coefficients of $\Phi$. E.g., if $d=3$, then $\I((\F,\F)_2)$ is the 
defining ideal of the twisted rational cubic curve. 

If $\II(a_0,\dots,a_d)$ is an invariant of degree $m$, then its 
{\sl evectant} is defined to be 
\begin{equation}
\E_\II = \frac{1}{m} \sum\limits_{i=0}^d \, 
\frac{\partial \II}{\partial a_i} \, (-x_2)^{d-i} \, x_1^i, 
\label{formula.evectant} \end{equation} 
which is a covariant of degree-order $(m-1,d)$. By Euler's formula we have an 
identity $(\E_\II,\F)_d= \II$. 

Let $\A \subseteq \rationals [a_0,\dots,a_d;x_1,x_2]$ 
denote the subring of covariants, which is naturally bigraded by $(m,q)$. 
By a fundamental theorem of Gordan, $\A$ is finitely generated. A 
minimal set of generators of $\A$ is called a fundamental system for $d$-ics. 
Moreover $\A$ is a unique factorization domain and 
each of the minimal generators is a prime element of $\A$.

The number of linearly independent covariants of $d$-ics of degree-order $(m,q)$ is given 
by the Cayley-Sylvester formula (see~\cite[Corollary 4.2.8]{Sturmfels}). 
For integers $n,k,l$, let $p(n,k,l)$ denote the number of partitions of 
$n$ into $k$ parts such that no part exceeds $l$. Then 
\begin{equation} \zeta_{m,q} = \dim \A_{m,q} = 
p \, (\frac{dm-q}{2},d,m)-p \, (\frac{dm-q-2}{2},d,m).  
\label{formula.CS} \end{equation} 

\begin{Example} \rm Let $d=5$, then 
$\zeta_{4,8} = p(6,5,4)-p(5,5,4) = 2$. 
A basis for the space $\A_{4,8}$ is given by 
\[ (\F,\F)_2 \, (\F,\F)_4, \quad \F \, (\F,(\F,\F)_4)_2. \] 
\end{Example} 

\subsection{Quintics} 
We will make use of the fundamental system for quintics, which has been known 
since the nineteenth century. The following table (adapted from~\cite[p.~131]{GrYo}) 
lists the degree-orders of the minimal generators of $\A$. For instance, there is 
one generator in degree-order $(5,3)$ and none in $(3,7)$. 

\[ \qquad \qquad \text{order} \] 
\[ \text{degree} \; \; 
\begin{tabular}{|c||c|c|c|c|c|c|c|c|c|} \hline 
{} &  0 &  1 &  2 &  3 & 4  &  5 &  6 &  7 &  9 \\ \hline \hline 
1  & {} & {} & {} & {} & {} &  1 & {} &  {} & {} \\ \hline 
2  & {} & {} &  1 & {} & {} & {} &  1 &  {} & {} \\ \hline 
3  & {} & {} & {} &  1 & {} & 1  & {} &  {} & 1 \\ \hline 
4  &  1 & {} & {} & {} &  1 & {} &  1 &  {} & {} \\ \hline 
5  & {} &  1 & {} &  1 & {} & {}  & {} &   1 & {} \\ \hline 
6  & {} & {} &  1 & {} &  1 & {} & {} &  {} & {} \\ \hline 
7  & {} &  1 & {} & {} & {} & 1  & {} &  {} & {} \\ \hline 
8  &  1 & {} &  1 & {} & {} &{}  & {} &  {} & {} \\ \hline 
9  & {} & {} & {} &  1 & {} &{}  & {} &  {} & {} \\ \hline 
11 & {} &  1 & {} & {} & {} &{}  & {} &  {} & {} \\ \hline 
12 &  1 & {} & {} & {} & {} &{}  & {} &  {} & {} \\ \hline 
13 & {} &  1 & {} & {} & {} &{}  & {} &  {} & {} \\ \hline 
18 & 1  & {} & {} & {} & {} &{}  & {} &  {} & {} \\ \hline 
\end{tabular} \] 

\vspace*{1cm} 
\medskip 

We will frequently need the following covariants: 
\begin{equation} 
\begin{array}{lll} 
\vartheta_{22}=(\F,\F)_4, & \vartheta_{26}=(\F,\F)_2, & 
\vartheta_{33}=(\vartheta_{22},\F)_2, \\ 
\vartheta_{39} =(\F,\vartheta_{26})_1, & 
\vartheta_{40}=(\vartheta_{22},\vartheta_{22})_2, & 
\vartheta_{44}=(\vartheta_{22},\vartheta_{26})_2, \\ 
\vartheta_{51}=(\vartheta_{22}^2,\F)_4, & 
\vartheta_{80}=(\vartheta_{22}^3,\vartheta_{26})_6.  
\end{array} \end{equation} 
The notation is so set up that $\vartheta_{m \,q}$ is a generator in 
degree-order $(m,q)$. (The comma is omitted for ease of reading.)

\medskip \medskip 

The computations which go into constructing such tables are generally very laborious, and of course 
the classical invariant theorists carried them out without the aid of machines. 
Hence, it is not unreasonable to worry about their correctness (also see the footnote 
on~\cite[p.~131-132]{GrYo}). In the case of binary quintics however, I have thoroughly checked 
that the table above is entirely correct. 

Here is a typical instance of how the table is used: we have 
\[ \zeta_{9,5} = p(20,5,9)-p(19,5,9) = 98-93 = 5, \] 
i.e., $\A_{9,5}$ is $5$-dimensional. Notice that 
\begin{equation} \label{basis.95}
B = \{ \vartheta_{51} \, \vartheta_{22}^2, \, 
\vartheta_{51} \, \vartheta_{44},  \, 
\vartheta_{40} \, \vartheta_{33} \, \vartheta_{22}, \, 
\vartheta_{40}^2 \, \F, \, \vartheta_{80} \, \F \} \end{equation}
are all of degree-order $(9,5)$. Since they are linearly independent over $\rationals$ (this can 
be checked by specializing to $F= x_1^5 + x_2^5 + (x_1+x_2)^5$ and solving a system of linear 
equations), $B$ is a basis of $\A_{9,5}$. This basis will be used in \S\ref{firstsyzygies.J}. 

Since $\zeta_{18,0} = p(45,5,18)-p(44,5,18) = 967-966=1$, up to scalar, quintics have 
a unique invariant of degree $18$. Hence, following~\cite[p.~131]{GrYo}, we will define 
\begin{equation} \HI = (\vartheta_{22}^7, \F \, \vartheta_{39})_{14}. 
\label{defn.H} \end{equation} 
(This merely requires checking that the transvectant is not identically zero, which can be 
done by specializing $\F$ and calculating directly.) Usually $\HI$ is called a skew-invariant 
(since it is of odd weight). Indeed, $\HI$ was the first discovery of a skew-invariant for any $d$. 
(They do not occur for $d \le 4$.) For what it is worth, a 
{\sc Maple} computation shows that $\HI$ is a linear combination of 
$848$ monomials in $a_0,\dots,a_5$. 

\subsection{} \label{geometry1} 
Let $u \in S_2$ be a nonzero vector. The duality in~(\ref{self-duality})
identifies the point $[u] \in \P S_2$ with its 
polar line $\{[v] \in \P^2: (u,v)_2 =0\} \in \P S_2^*$. 
The point lies on its own polar iff $(u,u)_2=0$, which happens iff $[u] \in C$. 
If $[u]$ lies on the polar of $[v]$, then $[v]$ lies on the 
polar of $[u]$. The pole of the line joining two points $[u],[v]$ is given by $[(u,v)_1]$. 
Three points $[u],[v],[w]$ are collinear iff $((u,v)_1,w)_2=0$. 

If $l \in S_1$, then the tangent to $\phi(l) \in C$ is the line 
$\{[l \, m]: m \in S_1 \}$. The line joining $\phi([l]),\phi([m])$ 
is (the polar of) $[l \, m]$. 

\subsection{} \label{R=H}
The following proposition will be needed in \S\ref{section.sing.H}. 
Let $G$ denote a binary quartic identified with four points 
$\Pi = \{\ssa,\ssb,\ssc,\ssd\} \subseteq C$. Consider the three pairwise 
intersections 
$\overline{\ssa \, \ssb} \cap \overline{\ssc \, \ssd}, 
\overline{\ssa \, \ssc} \cap \overline{\ssb \, \ssd}, 
\overline{\ssa \, \ssd} \cap \overline{\ssb \, \ssc}$, 
regarding each as a form in $S_2$. 
\begin{Proposition} \sl 
The product of the three points is given (of course up to scalar) by the 
covariant $\T(G) = (G,(G,G)_2)_1$. 
\label{prop.triple.intersection} \end{Proposition} 
\demo Let us write $G = a_\ux \, b_\ux \, c_\ux \, d_\ux$, where 
$a_\ux = a_1 \, x_1 + a_2\, x_2$ and $\ssa = \phi([a_\ux])$ etc. By \S\ref{geometry1}, the intersection 
$\overline{\ssa \ssb} \cap \overline{\ssc \ssd}$ corresponds to 
\[ ((a_\ux^2,b_\ux^2)_1,(c_\ux^2,d_\ux^2)_1)_1 = 
(a \, b)(c \, d) \, (a_\ux \, b_\ux, c_\ux \, d_\ux)_1,  \] 
where $(a \, b) = a_1 b_2 - a_2 \, b_1$ etc. Hence, up to a factor, the product corresponds to 
\begin{equation} 
(a_\ux \, b_\ux, c_\ux \, d_\ux)_1 \, (a_\ux \, c_\ux, b_\ux \, d_\ux)_1 \,  
(a_\ux \, d_\ux, b_\ux \, c_\ux)_1.  
\label{product.3} \end{equation} 
The last expression is of degree $3$ in the coefficients of $G$ (since each of the letters 
$a, \dots, d$ occurs thrice), moreover it is a covariant since the underlying 
geometric construction is compatible with the $SL(V)$-action. However, $\zeta_{3,6} = 1$ for 
binary quartics, hence $\T(G)$ and (\ref{product.3}) are equal up to a scalar. \qed 

\smallskip 

The result remains true if $\Pi$ contains one double point, say $\ssa = \ssb$, with
$\overline{\ssa \, \ssb}$ interpreted as the tangent to $C$ at $\ssa$. 
By~\cite[\S 3.5.2]{Glenn}, the covariant $\T(G)$ vanishes identically iff $\Pi$ consists of two (possibly coincident) double points, 
say $\ssa = \ssb, \ssc = \ssd$. In this case the 
geometric construction collapses, since $\overline{\ssa \, \ssc} 
\cap \overline{\ssb \, \ssd}$ is no longer a determinate point. 

\smallskip 

This proposition can be used to give an alternate definition of $\HI$. 
Let $\R$ denote the resultant $\text{Res}(\F,\vartheta_{33})$, defined 
as the determinant of an $8 \times 8$ Sylvester matrix 
(see~\cite[Ch.~V,\S 10]{Lang}). 
By construction it is of degree $5 \times 3 + 3 \times 1 =18$ in the $\{a_i\}$. 

\begin{Proposition} \sl 
The hypersurface defined by $\R$ coincides with $\He$. 
\end{Proposition} 
We will avoid using the fact that $\zeta_{18,0} =1$. 

\smallskip 

\demo 
Let us first show that $\R$ is not identically zero. Specialize to 
$F = x_1^5 + 2 \, x_2^5 + (x_1+ x_2)^5$. Then 
$\vartheta_{33} = -12 \, x_1 \, x_2 \, (x_1+x_2)$, 
which has no common factor with $F$, hence $\R \not\equiv 0$. 
Now assume that $F$ and $\vartheta_{33}(F)$ have a 
common linear factor, we may take it to be $x_1$ after a change of variables. 
Let $F = x_1 \, G$, with $G = (a_0,a_1,a_2,a_3,a_4 \cb x_1,x_2)^4$. 
Calculating directly, we have 
\begin{equation} 
\vartheta_{33}(F)|_{x_1:=0} = 
\frac{24}{125} \, x_2^3 \, (2 \, a_3^3 + a_1 \, a_4^2 - 3 \, a_2 \, a_3 \, a_4), 
\label{vartheta33.1} \end{equation} 
which vanishes by hypothesis. Hence 
\[ \T(G)|_{x_1:=0}  = 
- \, x_2^6 \, (2 \, a_3^3 + a_1 \, a_4^2 - 3 \, a_2 \, a_3 \, a_4) \] 
must also vanish, i.e., $x_1$ must divide one of the three 
intersection points coming from $G$. Denote this point by $\ssz = 
[x_1 \, (\alpha \, x_1 + \beta \, x_2)]$. It is now immediate that 
the divisor corresponding to $F$ is in involution with respect to $\ssz$, hence 
$[F] \in \He$. Thus we have an inclusion of hypersurfaces 
$\{[F] \in \P^5: \R = 0 \} \subseteq \He$. Since the latter 
is irreducible, they must be equal. \qed 

\subsection{} 
It will prove useful to introduce the following loci in $\P^5$. 
If $\lambda=(\lambda_1,\dots,\lambda_r)$ is a partition of $5$, 
let $X_\lambda$ denote the closed subvariety 
\[ \{[F] \in \P^5: F = \prod l_i^{\lambda_i} \; \; 
\text{for some $l_i \in S_1$} \}. \] 
In other words, the divisor of $[F] \in X_\lambda$ is of the form 
$\lambda_1 \ssa_1 + \dots + \lambda_r \, \ssa_r$ with some of the 
$\ssa_i$ possibly coincident. The dimension of $X_\lambda$ equals the 
number of (nonzero) parts in $\lambda$. There is an inclusion 
$X_\mu \subseteq X_\lambda$ iff $\lambda$ is a refinement of $\mu$. 
For instance, $X_{(5)}$ is the rational normal 
quintic, $X_{(2,1,1,1)}$ is the discriminant hypersurface, and 
$X_{(3,1,1)}$ is the locus of nullforms. 

\subsection{A summary of results} \label{section.results.summary}
In \S\ref{section.sing.H} we will construct a desingularization of $\He$, and then show that 
its singular locus $\B$ consists of three components 
$\Omega_{(1)},\Omega_{(2)}$ and $X_{(3,1,1)}$. They are respectively the 
$SL_2$-orbit closures of the forms 
\[ x_1^5 + x_2^5, \quad 
x_1 \, x_2 \, (x_1-x_2) \, (x_1^2 + x_1 \, x_2 + x_2^2), \quad 
x_1^3 \, x_2 \, (x_1+x_2). \]
Their degrees are $6,10$ and $9$, hence $\B$ is of degree $25$ 
and pure codimension two. Next we show that the ideal 
$I_\B \subseteq R$ is a complete intersection, defined by the coefficients 
of $\vartheta_{51}$. 

In \S\ref{section.dual.H} it will be seen that $\He$ is naturally isomorphic to its own 
dual variety. The duality $S_5 \simeq S_5^*$ in (\ref{self-duality})
induces an isomorphism $\sigma: \P^5 \stackrel{\sim}{\lra} (\P^5)^*$. 
Let $[F] \in \He \setminus \B$, with $T_{\He,[F]}$ the tangent space to 
$\He$ at $[F]$. Then the point $\sigma^{-1}(T_{\He,[F]})$ coincides with 
$[\E_\HI(F)]$ (the value of the evectant at $F$). It turns out however, 
that this point also belongs to $\He$. Thus we get a morphism 
\[ \He \setminus \B \lra \He \setminus \B, \quad 
[F] \lra [\E_\HI(F)]. \] 
This map is involutive, i.e., $\E_\HI(\E_\HI(F))$ equals $F$ up to a scalar. 

Let $J = (\frac{\partial \, \HI}{\partial a_0}, \dots, 
\frac{\partial \, \HI}{\partial a_5}) \subseteq R $
denote the Jacobian ideal of $\HI$. In~\S\ref{section.jacobian} 
we show that $J$ is a perfect ideal of height two, with an $SL_2$-equivariant 
minimal resolution 
\[ \begin{aligned} 
0 \la R/J \la R & \la R(-17) \otimes S_5 \\ 
& \la R(-18) \otimes S_2 \oplus R(-22) \oplus R(-26) \la 0. 
\end{aligned} \] 
During the course of the proof we will see that $J$ naturally fits into 
a three-parameter family of perfect ideals. 

The results of~\S\ref{section.J} allow us to identify the morphisms in this resolution 
up to three distinct possibilities, but no further. 
In order to resolve this ambiguity it would suffice to 
calculate the value of $\E_\HI$ at $\cF_Q$. A general formalism is developed in 
\S\ref{section.evectants} to solve this problem. For any covariant 
$\Phi$ of $d$-ics, we construct a sequence of covariants $\A_\bullet$ called its evectants; this 
generalizes the classical construction from \S\ref{section.evectant}. Given two arbitrary 
covariants $\Phi,\Psi$ with evectants $\A_\bullet,\B_\bullet$, we deduce formulae for 
calculating the evectants of a general transvectant $(\Phi,\Psi)_r$. 
This iterative scheme is then applied to formula~(\ref{defn.H}) to evaluate $\E_\HI$. 
Nearly all of \S\ref{section.evectants} can be read independently of the rest of the paper. 
\label{end.summary}

\subsection{A note on computational procedures} Since I have used 
machine computations in several parts of this paper, their role 
and extent should be clarified. All the computations have been done in 
{\sc Maple}. I have written routines to calculate the numbers $p(n,k,l)$ and 
$\zeta_{m,q}$ appearing in formula (\ref{formula.CS}). I have also programmed 
formula~(\ref{trans.formula}) for calculating transvectants; hence 
identities such as~(\ref{expr.TA}) and (\ref{phi51.1}) are machine-computed. 
I have also used {\sc Maple} for some routine calculation in linear 
algebra, e.g., for evaluating Wronskian determinants and for solving systems of linear equations. 
None of the results depend upon calculating Gr{\"o}bner bases in 
any guise (e.g., minimal free resolutions). 

\setcounter{footnote}{1}
On the whole, I have not succeeded in bypassing heavy calculations entirely, 
and I very much doubt if this is at all possible. The Hermite invariant is a specific 
algebro-geometric object which is not a member of any natural `family', hence it seems 
unlikely that merely general considerations will enable us to prove much about it. Even so, 
I believe that none of the calculations done here by a machine are beyond the ambit of a 
patient and able human mathematician\footnote{Paul Gordan and George Salmon come to mind; for 
instances, see~\cite{Gordan} or the tables at the end of~\cite{Salmon1}.}.

\section{The singular locus} \label{section.sing.H}
\subsection{} First we construct a natural desingularization of $\He$. 
Let 
\[ Y = \{ (\ssc,\ssz) \in C \times \P^2: 
\text{the tangent to $C$ at $\ssc$ passes through $\ssz$} \}. \] 
The second projection $Y \stackrel{\alpha}{\lra} \P^2$ is a double 
cover ramified along $C$. Let 
$\P \, T_{\P^2} \lra \P^2$ denote the projectivisation of the tangent 
bundle of $\P^2$, so that the fibre over $\ssz \in \P^2$ can be identified 
with the pencil of lines through $\ssz$. Define the $\P^2$-bundle 
\[ \text{Sym}^2 \, (\P \, T_{\P^2}) 
\stackrel{\beta}{\lra} \P^2, \] 
so that an element in $\beta^{-1}(\ssz)$ is an unordered pair of 
(possibly coincident) lines $L_1,L_2$ through $\ssz$. 
Consider the pullback square 
\[ \diagram 
{\cZ \,} \rto \dto & \text{Sym}^2 \, (\P \, T_{\P^2}) \dto^\beta \\ 
Y \rto_\alpha & {\P^2  \, .}
\enddiagram \] 
Define $\cZ \stackrel{f}{\lra} \He$ 
by sending $(\ssc,\ssz) \times (L_1,L_2)$ to the divisor 
\[ \ssc + L_1 \cap C + L_2 \cap C. \] 
(Of course, $L_i \cap C$ are interpreted scheme-theoretically.) 
By construction $f$ is a projective birational morphism which is a 
desingularization of $\He$. We will use this map to detect the 
singularities of $\He$. 
Since $Y$ is a rational variety (in fact isomorphic to $\P^1 \times \P^1$), 
so is $\cZ$ and hence $\He$. 
Henceforth we will write 
$(\ssc,\ssz;L_1,L_2)$ for $(\ssc,\ssz) \times (L_1,L_2) \in \cZ$. 

\begin{Lemma} \sl 
The morphism 
$\cZ \setminus f^{-1}(X_{(5)}) \stackrel{f}{\lra} 
\He \setminus X_{(5)}$ is finite. 
\end{Lemma} 
\demo Since the morphism is projective, it suffices to show that it has 
finite fibres (see~\cite[Lemma 14.8]{Harris}). 
Let $(\ssc,\ssz;L_1,L_2) \in f^{-1}([F])$. There are finitely many 
choices for $\ssc$. By hypothesis there is a point 
$\ssa (\, \neq \ssc)$ appearing in $[F]$; hence for a given $\ssc$ 
there are only finitely many possibilities for $\ssz$ (because 
$\overline{\ssz \, \ssa} \cap C$ must be contained in $[F]$). Then for a 
given $\ssz$, there are only finitely many possibilities for the $L_i$. 
\qed 

\smallskip 

This argument breaks down over $X_{(5)}$; in fact 
$f^{-1}(X_{(5)}) \lra X_{(5)}$ is a $\P^1$-bundle. 

\subsection{} 
Define the forms 
\[ \begin{array}{ll} 
\bU_{(1)} = x_1^5 + x_2^5, & \bU_{(2)} = x_1 \, x_2 \, (x_1 -x_2) \, (x_1^2 + x_1 \, x_2 + x_2^2)  \\ 
\bU_{(3)} = x_1^3 \, x_2 \, (x_1 + x_2), & \bU_{(4)} = x_1^3 \, x_2^2, \\ 
\bU_{(5)} = x_1^4 \, x_2, & \bU_{(6)} = x_1^5. 
\end{array} \] 
Let $\B \subseteq \He$ denote the union of the orbits of all the $U_{(i)}$. 
We claim that $\B$ is closed. Indeed, by~\cite[\S2]{Aluffi-Faber} the closure of any orbit is a union of 
orbits of forms of the type $x_1^a \, x_2^b$, and they are already included. 

\begin{Theorem} \sl 
The singular locus $\text{Sing}(\He)$ coincides with $\B$. 
\label{theorem.sing} \end{Theorem} 

The theorem will follow from the following proposition. 
\begin{Proposition} \sl 
\begin{enumerate} 
\item 
For $[F] \in \He$, the fibre $f^{-1}([F])$ consists of more than one point iff 
$F$ lies in the orbit of one of the forms $U_{(i)}$ for $1 \le i \le 6,i \neq 5$. 
\item 
Assume $[F] \in \He \setminus \B$, and $f^{-1}([F]) = \{\ssw\}$. Then the morphism on tangent spaces 
$T_{\cZ,\ssw} \lra T_{\He,[F]}$ is injective. 
\end{enumerate} 
\end{Proposition} 

Let us show the theorem assuming the proposition. If $[F]$ lies in the orbit of one of 
$\bU_{(1)},\dots,\bU_{(4)}$, then the fibre $f^{-1}([F])$ is disconnected, hence 
$[F]$ is not a normal point. Since $\bU_{(5)},\bU_{(6)}$ lie in the orbit closure of 
$\bU_{(3)}$, we deduce that $\B \subseteq \text{Sing}(\He)$. If $[F] \in \He \setminus \B$, then 
by~\cite[Theorem 14.9]{Harris} the map $f$ is a local isomorphism in a neighbourhood of $\ssw$, 
hence $[F]$ is a nonsingular point. \qed 

\subsection{} Let us prove part (1) of the proposition. Define 
\[ \cS = \{[F] \in \He: f^{-1}([F]) \; \text{consists of at least two points} \}. \]
Evidently $\bU_{(6)} \in \cS$. 
Assume that $[F] = 3 \ssc + \ssa_1 + \ssa_2$, where $\ssa_1,\ssa_2$ are (possibly coincident) points each 
different from $\ssc$. Let $\ssz$ denote the intersection $\overline{\ssc \ssc} \cap 
\overline{\ssa_1 \ssa_2}$, then 
$(\ssc,\ssz; \overline{\ssc \, \ssc},\overline{\ssa_1 \, \ssa_2})$ and 
$(\ssc,\ssc; \overline{\ssc \, \ssa_1},\overline{\ssc \, \ssa_2})$ both map to 
$[F]$; this shows that $\bU_{(3)}, \bU_{(4)} \in \cS$. It is equally clear that 
$\bU_{(5)} \notin \cS$. 

If a point of the form $(\ssc,\ssc,L_1,L_2)$ belongs to $f^{-1}([F])$, then $[F]$ must have a 
point of multiplicity $\ge 3$ at $\ssc$, which is already considered above. 
Hence assume that $[F] \in \cS \setminus X_{(3,1,1)}$, and 
$(\ssc,\ssz;L_1,L_2),(\ssc',\ssz';L_1',L_2')$ are two distinct points in $f^{-1}([F])$. 
Since $\ssc \neq \ssz$, we may write $\ssc = \phi([x_1]), \ssz = [x_1 x_2]$ after a 
change of variables. Then $[F] = [\cF_Q]$ for some $Q \in \P^2$ (see~\S\ref{section.FQ}). 

If $q_0=0$, then both $q_1,q_2$ must be nonzero (otherwise $[\cF_Q] \in X_{(3,1,1)}$). 
But then $[\cF_Q]$ is in the orbit of $A=x_1 \, x_2^2 \, (x_1+x_2) \, (x_1-x_2)$, and it is 
clear from the geometry that $[A] \notin \cS$. 

Hence we may assume $q_0 =1$, and then 
\[ F = x_1 \, (x_1 - \alpha \, x_2) \, (x_1 + \alpha \, x_2) \, 
(x_1 - \beta \, x_2) \, (x_1 + \beta \, x_2) \] 
for some $\alpha,\beta$, such that $\ssc' = \phi([x_1 - \alpha \, x_2])$.  
By assumption $\ssz'$ is one of the diagonal intersection points (see~\S\ref{R=H}) coming from 
the quartic form $G = x_1 \, (x_1 + \alpha \, x_2) \, 
(x_1 - \beta \, x_2) \, (x_1 + \beta \, x_2)$. The quadratic form corresponding to 
$\ssz'$ must divide $\T(G)$, and hence $x_1 - \alpha \, x_2$ must divide $\T(G)$. 
By a direct calculation, 
\begin{equation} \begin{aligned} 
{} & \T(G)|_{x_1 := \alpha \, x_2} \\ 
= \; & \frac{1}{32} \, x_2^6 \, \alpha^3 \, (\alpha^2 + 3 \, \beta^2) \, 
(\alpha^2 + 4 \, \alpha \, \beta - \beta^2) \, 
(\alpha^2 - 4 \, \alpha \, \beta - \beta^2), 
\end{aligned} \label{expr.TA} \end{equation}
which must vanish. Now $\alpha \neq 0$, since $[F] \notin X_{(3,1,1)}$. 
Hence we have two cases 
\[ \frac{q_0 \, q_2}{q_1^2} = 
\frac{4 \, \alpha^2 \, \beta^2}{(\alpha^2 + \beta^2)^2}  = 
\begin{cases} 
1/5 & \text{if $\alpha^2 \pm 4 \, \alpha \beta - \beta^2 =0$,} \\ 
-3 & \text{if $\alpha^2 + 3 \, \beta^2 =0$.}
\end{cases} \] 
A form satisfying the first case is in the orbit of 
$\cF_{[1,5,5]} = x_1 (1,5,5 \cb x_1^2,x_2^2)^2$. By the transformation 
$(x_1,x_2) \lra (x_1 +x_2,x_1-x_2)$ it can be brought into the 
more manageable form 
\begin{equation} \bU_{(1)} = x_1^5 + x_2^5. 
\label{psi.1} \end{equation} 
Similarly in the second case $\cF_{[1,1,-3]}$ can be brought into the form 
\begin{equation} 
\bU_{(2)} = x_1 \, x_2 \, (x_1-x_2) \, (x_1^2 + x_1 \, x_2 + x_2^2). 
\label{psi.2} \end{equation}
via $(x_1,x_2) \lra (x_1-x_2,x_1+x_2)$. We have shown that any form in 
$\cS \setminus X_{(3,1,1)}$ belongs to the orbit of either $\bU_{(1)}$ or $\bU_{(2)}$. 
It remains to show that the latter two belong to $\cS$, this 
can be done by an explicit construction as follows: 

Let $\omega = \exp(\frac{\pi \sqrt{-1}}{5})$, and $y = \omega^r \, x_2$. 
Define points $\ssc = \phi([x_1-y]), \ssz = [(x_1+y) \, (x_1-y)]$, 
and $L_i$ to be the line joining 
$\phi([x_1-\omega^i \, y])$ and $\phi([\omega^i \, x_1-y])$ 
for $i =1,2$. This 
gives a point of $f^{-1}([\bU_{(1)}])$ for every $1 \le r \le 5$. 

Let $\nu = \exp(\frac{2 \pi \sqrt{-1}}{3})$, and $y = \nu^r \, x_2$. 
Define points $\ssc = \phi([x_1-y]), \ssz = [(x_1-y) \, (x_1+y)]$. 
Let $L_1$ be the line joining $\phi([x_1]),\phi([x_2])$, and 
$L_2$ joining $\phi([x_1-\nu \, y])$ and $\phi([\nu \, x_1 - y])$. 
This gives a point of $f^{-1}([\bU_{(2)}])$ for every $1 \le r \le 3$. 

This completes the proof of part (1). \qed 

\subsection{} We will prove part (2) by introducing a local parametrisation of the affine 
version of $f$,  and directly calculating the map on tangent spaces. Since 
$[F] \notin X_{(3,1,1)}$, after a change of variables we may write $F = 
x_1 \, (1,\xi,1 \cb x_1^2,x_2^2)^2$ for some $\xi \in \complex$. 

Let $\A = S_1 \times S_1 \times \complex$, and define a morphism 
from $\A$ to $\cZ$ by sending $(l_1,l_2,\xi) \in \A$ 
to $([l_1^2],[l_1 \, l_2],L_1,L_2)$, where $L_1, \, L_2$ correspond to the solutions of 
the equation $(1,\xi,1 \cb l_1^2,l_2^2)^2=0$. Since the morphism is smooth, for a local 
parametrisation of $f$ we may use the map 
\[ \hf: \A \lra \text{Cone}(\He), \quad 
(l_1,l_2,\xi) = l_1 \, (1,\xi,1 \cb l_1^2,l_2^2)^2.  \] 
The image of an arbitrary tangent vector $(m_1,m_2,\eta)$ via $d \hf$ 
is given by the limit 
\[ \begin{aligned} 
{} & \tau(m_1,m_2,\eta) \\ 
= & \lim_{\epsilon \ra 0} \; 
\frac{1}{\epsilon} \, [ \, 
\hf(l_1 + \epsilon \, m_1,l_2 + \epsilon \, m_2,\xi + \epsilon \, \eta) - 
\hf(l_1,l_2,\xi) \, ].  
\end{aligned} \]  

Writing $\ssw = (x_1,x_2,\xi)$, the image of the map $T_{\A,\ssw} \lra T_{\text{Cone}(\He),F}$ is 
spanned by the five vectors 
\[ \begin{array}{ll} 
\tau(x_1,0,0) = x_1 \, (5,3 \, \xi, 1 \cb x_1^2, x_2^2), & 
\tau(x_2,0,0) = x_2 \, (5,3 \, \xi, 1 \cb x_1^2, x_2^2), \\ 
\tau(0,x_1,0) = 4 \, x_1^2 \, x_2 \, (\xi \, x_1^2 + x_2^2), & 
\tau(0,x_2,0) = 4 \, x_1 \, x_2^2 \, (\xi \, x_1^2 + x_2^2), \\ 
\tau(0,0,1) = 2 \, x_1^3 \, x_2^2. 
\end{array} \] 
In order to verify that they are linearly independent, we calculate their Wronskian 
\begin{equation} \begin{aligned} 
{} & \left| \begin{array}{rrrrr} 
600\, x_1 & 72\, \xi \, x_2 &  72\, \xi x_1\, & 24\, x_2 & 24\, x_1\\ 
120\, x_2 & 120\, x_1 & 72\, \xi\, x_2 & 72\, \xi \, x_1 & 120\, x_2 \\
96\, \xi\, x_2 & 96\, \xi \, x_1 & 48\, x_2 & 48\, x_1 & 0 \\ 
0 & 48\, \xi\, x_2 & 48\, \xi \, x_1 & 96\, x_2 & 96\, x_1 \\ 
0 & 24\, x_2 & 24\, x_1 & 0 & 0 
\end{array} \right| \\ 
= & - \, 2^{18} \, 3^5 \, 5^2 \, x_1 \, 
(6 \, \xi^2 -5, -5 \, \xi, 5 \cb x_1^2,x_2^2)^2. 
\end{aligned} \label{formula.wr=E} \end{equation} 
This is nonzero for any $\xi$, which proves part (2) 
of the proposition. The proof of 
Theorem~\ref{theorem.sing} is complete. \qed 

\medskip 

One can restate the theorem as follows: 
$\cF_Q$ is a singular point of $\He$, iff one of the expressions 
$q_2, q_0 \, q_2 + 3 \, q_1^2,5 \, q_0 \, q_2 - q_1^2$ is zero. 

\subsection{} For $i=1,2$, let $\Omega_{(i)}$ denote the orbit closure of $[\bU_{(i)}]$, and 
let $\G_i \subseteq SL(V)$ denote the stabilizer subgroup of 
$[\bU_{(i)}]$. By~\cite[\S0]{Aluffi-Faber}, we have a formula 
\[ \deg \, \Omega_{(i)} = \frac{5.4.3}{|\G_i|}. \] 
Since an element of $\G_i$ must permute the linear factors of $\bU_{(i)}$, it is 
easy to determine all symmetries by mere inspection. 
The group $\G_1$ is the dihedral group $D_5$ of order $10$, generated by 
the transformations 
\[ (x_1,x_2) \lra 
\begin{cases} (x_2,x_1), \\  (x_1,\exp(\frac{2 \pi \sqrt{-1}}{5} ) \, x_2). 
\end{cases} \] 
Similarly $\G_2$ is isomorphic to $D_3$, generated by 
\[ (x_1,x_2) \lra 
\begin{cases} (x_2,x_1), \\ (\exp(\frac{2 \pi \sqrt{-1}}{3}) \, x_1,x_2). \end{cases} \] 
Hence $\Omega_{(1)},\Omega_{(2)}$ are of degrees $6$ and $10$ 
respectively. The degree of $X_{(3,1,1)}$ is $9$, as 
given by a formula due to Hilbert~\cite{Hilbert1}. 
\subsection{} 
Let $\fp_{(i)} \subseteq R$ denote the homogeneous ideal of 
$\Omega_{(i)}$. The variety $\Omega_{(1)}$ is the closure of the 
union of secant lines to $X_{(5)}$, and it is known 
(as an instance of a more general result) that $\fp_{(1)}$ is a perfect ideal 
of height two (see~\cite[Theorem~1.56]{Iarrobino-Kanev}). We briefly recapitulate the proof. 
Given $F \in S_5$, define 
\[ \alpha_F: S_2 \lra S_3, \quad G \lra (F,G)_2, \] 
and let 
\[ \alpha: S_2 \otimes R(-1) \lra S_3 \otimes R \] 
denote the corresponding morphism of graded $R$-modules (\S\ref{section.map.Rmodules}). 

\begin{Lemma} \sl 
The map $\alpha_F$ is injective for a general $F$, moreover $\ker \alpha_F$ is nonzero 
iff $[F] \in  \Omega_{(1)}$. 
\end{Lemma} 
\demo It is easily verified from formula~(\ref{trans.formula}) that 
$\ker \alpha_F = 0$ for $F = x_1^5 + x_2^5 + (x_1+x_2)^5$. Assume $G (\neq 0) \in \ker \alpha_F$, 
then after a change of variables $G$ can be written as either $x_1^2$ or $x_1 \, x_2$. 
In the former case $F = x_1^4 \, (c_1 \, x_1 + c_2 \, x_2)$ and in the latter case $F =  c_1 \, x_1^5 + 
c_2 \, x_2^5$. The `if' part is equally clear. \qed 

\smallskip 

By the Porteous formula (see~\cite[Ch.~II.4]{ACGH}) the scheme-theoretic degeneracy locus 
$\{\rank \alpha_F \le 2 \}$ has degree $6$ (it is the coefficient of $h^2$ in the Maclaurin 
expansion of $(1+h)^{-3}$), and so does $\Omega_{(1)}$. Hence 
the ideal of maximal minors of $\alpha$ coincides with $\fp_{(1)}$, and we get a 
Hilbert-Burch resolution (see~\cite[\S20.4]{Ei})
\[ 0 \la R/\fp_{(1)} \la R \stackrel{\delta_0}{\la} R(-3) \otimes S_3 
\stackrel{\delta_1}{\la} R(-4) \otimes S_2 \la 0.  \] 
Now consider the complex 
\[ R \stackrel{\delta_0^\vee}{\ra} R(3) \otimes S_3 
\stackrel{\delta_1^\vee}{\ra} R(4) \otimes S_2.  \] 
To describe the first map, let $\W_{(1)}$ denote the Wronskian of $\alpha_\F$, 
i.e., the determinant of the $3 \times 3$ matrix of linear forms 
\[ (i,j) \lra \frac{\partial^2 \,  (\F,x_1^{3-i} \, x_2^{i-1})_2}{\partial x_1^{3-j} \, \partial x_2^{j-1}}, \quad 
(1 \le i,j \le 3). \] 
Now $\W_{(1)}$ is a covariant of degree-order $(3,3)$, and $\zeta_{3,3}=1$ for quintics, 
hence it must coincide with $\vartheta_{33}$ up to a scalar. Thus $\fp_{(1)} = \I(\vartheta_{33})$. 

Up to a scalar, the map $\delta_1^\vee$ must be given by 
$S_3 \lra S_2, G \lra (F,G)_3$. From $\delta_1^\vee \circ \delta_0^\vee=0$ 
we deduce the identity $(\vartheta_{33},\F)_3=0$. 

\subsection{} Using similar ideas we will find a free resolution of $\fp_{(2)}$. It is sensible to look for 
a $4 \times 5$ matrix of linear forms, since then by Porteous's formula 
the degeneracy locus $\{ \rank \le 3\}$ has expected degree $10$. 
\begin{Proposition} \sl The ideal $\fp_{(2)}$ is perfect of height two. 
\end{Proposition} 
\demo Consider the map 
\[ \beta_F: S_3 \lra S_4, \quad G \lra (F,G)_2, \] 
and let $\W_{(2)}$ denote the corresponding $4 \times 4$ Wronskian determinant 
\[ (i,j) \lra \frac{\partial^3 \, (\F,x_1^{4-i} \, x_2^{i-1})_2}
{\partial x_1^{4-j} \, \partial x_2^{j-1}}, \qquad 
(1 \le i,j \le 4), \] 
which is a covariant of degree-order $(4,4)$. 
Let ${\mathfrak a} = \I(\W_{(2)})$ denote the ideal of maximal minors; {\sl a priori} we know it 
to be of height $\le 2$. If it were to have height one, then an invariant would have to divide $\W_{(2)}$, 
which is impossible. Hence we get a free resolution 
\[ 0 \la R/{\mathfrak a} \la R \la R(-4) \otimes S_4 \la R(-5) \otimes S_3 \la 0. \] 
Now a direct calculation shows that 
\[ \begin{aligned} 
\W_{(2)}(\cF_Q) & = \left| \begin{array}{rrrr}
24/5 \, q_1 \, x_1 & 12/5 \, q_2 \, x_2 & 12/5 \, q_2 \, x_1 & 0 \\ 
-2 \, q_1 \, x_2 & -2 \, q_1 \, x_1 & 2/5 \, q_2 \, x_2 & 
2/5 \, q_2 \, x_1\\ 
8 \, q_0 \, x_1 & -4/5 \, q_1 \, x_2 & -4/5 \, q_1 \, x_1 & 
-16/5 \, q_2 \, x_2 \\ 
6 \, q_0 \, x_2 & 6 \, q_0 \, x_1 & 18/5 \, q_1 \, x_2 & 18/5 \, q_1 \, x_1 
\end{array} \right| \\ 
& = \frac{1152}{125} \, (q_0 \, q_2 + 3 \, q_1^2) \, 
(5 \, q_0 \, q_2 + q_1^2, -2 \, q_1 \, q_2, 2 \, q_2^2 \cb x_1^2,x_2^2)^2. 
\end{aligned} \] 
Hence $\W_{(2)}$ vanishes on $\Omega_{(2)}$. Since the latter has degree $10$, the 
scheme defined by ${\mathfrak a}$ coincides with $\Omega_{(2)}$ and ${\mathfrak p}_{(2)} = 
{\mathfrak a}$. \qed 

\smallskip 

A basis for the space $\A_{4,4}$ is given by the two covariants $\vartheta_{22}^2,
\vartheta_{44}$, hence $\W_{(2)}$ must be their linear combination. 
The actual coefficients can be easily found by specializing $F$ and 
then solving a system of linear equations. This gives the relation 
$\W_{(2)} = 1/5760 \, (7 \, \vartheta_{22}^2 - 10 \, \vartheta_{44})$. 
As before, we have an identity $(\W_{(2)},\F)_3=0$. 

\subsection{} 
By a result of Weyman (see~\cite[Theorem 3]{Weyman}), 
the ideal of $X_{(3,1,1)}$ 
(say $\fq$) is generated in degrees $\le 4$. If we specialize to 
$F = x_1^3 \, x_2 \, (x_1+x_2)$ and search through 
all covariants in degrees $\le 4$, then we find that only 
$\vartheta_{40}$ and $2 \, \vartheta_{22}^2 + 15 \, \vartheta_{44}$ 
vanish on $F$, hence their coefficients must generate $\fq$. One sees
that $\fq$ is not perfect; indeed, it would have to arise as the ideal of 
maximal minors of a map \[ R \otimes (S_0 \oplus S_4) \lra 
\bigoplus\limits_{i \ge 0}  \, R(i) \otimes 
(S_{k_i} \oplus S_{k_i'} \oplus \dots)  \] 
such that the target module has rank $5$, the minors are of degree $4$ and the Porteous 
degree is $9$. However no such integers can be found. 

\subsection{} Let $I_\B \subseteq R$ denote the defining ideal of the singular locus $\B$. 
\begin{Proposition} \sl The ideal $I_\B$ is a complete intersection generated by the 
two coefficients of the covariant $\vartheta_{51}$. 
\end{Proposition} 
\demo The ideal ${\mathfrak e} = \I(\vartheta_{51})$ is a complete intersection, since otherwise an 
invariant would have to divide both coefficients of $\vartheta_{51}$. By a direct calculation, 
\begin{equation} \vartheta_{51}(\cF_Q) = 
\frac{4}{625} \, q_2 \, (q_0 \, q_2 + 3 \, q_1^2) \, (5 \, q_0 \, q_2 - q_1^2) \, x_1, 
\label{phi51.1} \end{equation} 
hence $\vartheta_{51}(F)$ vanishes on $\B$. Since $\deg \B = 25$, we must have 
${\mathfrak e} = I_\B$. \qed 

\smallskip 

Given a point $[F] \in \He \setminus \B$, the linear form $\vartheta_{51}$ `detects' the point of 
tangency $\ssc$ in the configuration on page~\pageref{conic_diagram}. 
Indeed this is visibly true of $\cF_Q$, and since $\vartheta_{51}$ is a covariant, it is true generally. 

\section{The dual variety} \label{section.dual.H}
Let $\sigma: \P S_d \stackrel{\sim}{\lra} \P S_d^*$ be the isomorphism 
induced by the duality  in (\ref{self-duality}); it identifies $[A] \in \P S_d$ with the 
hyperplane $\{[B] \in \P S_d: (A,B)_d =0 \}$. 
Let $\II$ denote a degree $m$ invariant of $d$-ics, defining a hypersurface 
$\X \subseteq \P S_d$. 

\begin{Proposition} \sl 
Let $[F] \in \X$ be a nonsingular point, and let $T = T_{\X,[F]} \in \P S_d^*$ 
denote the tangent space to $\X$ at $[F]$. Then we have an equality 
\[ [\E_\II(F)] = \sigma^{-1}(T). \] 
\end{Proposition} 
\demo Let $B = (b_0,\dots,b_d \cb x_1,x_2)^d$. The point 
$[b_0,\dots,b_d]$ belongs to $T$ iff 
\[ \sum\limits_{i=0}^d \, b_i \, 
(\left. \frac{\partial \II}{\partial a_i}\right|_F) = 0. \] 
This condition can be rewritten as $(\E_\II(F),B)_d =0$, hence the assertion. \qed

Now let $F = x_1 \, (1,\xi,1 \cb x_1^2,x_2^2)^2$. By the Proposition together with 
Lemma~\ref{lemma.wr}, the evectant $\E_\HI(F)$ is 
given (up to scalar) by the Wronskian of a basis of $T_{\He,[F]}$. But 
we have already calculated the latter in~(\ref{formula.wr=E}). 
After the substitution 
\[ (x_1,x_2,\xi) \lra (q_0^{1/5} \, x_1, \, q_2^{1/4} \, q_0^{-1/20} \, x_2, 
\, q_1 \, q_0^{-1/2} \, q_2^{-1/2}) \] 
we get the expression 
\[ \E_\HI(\cF_Q) = \text{constant} \times \cF_{Q'}, \] 
where 
\[ Q' = [\, q_0 \, q_2 - \frac{6}{5} \, q_1^2, q_1 \, q_2, -q_2^2 \, ]. \] 
Since $\E$ is a degree $17$ covariant, the `constant' must be a degree 
$15$ polynomial in the $q_i$. Now $\E_\HI(F)$ vanishes identically iff $[F] \in \B$, 
so we must have 
\begin{equation} 
\E_\HI(\cF_Q) = \constant \, q_2^n \, (q_0 \, q_2 + 3 \, q_1^2)^{n'}
\, (5 \, q_0 \, q_2 - q_1^2)^{n''} \cF_{Q'}, 
\label{EHI.exp1} \end{equation} 
for some integers $n,n',n''$ such that $n + 2 \, n' + 2 \, n''=15$. Here (and subsequently) 
$\constant$ stands for some {\sl nonzero} rational number which need not be precisely specified. 
The indices $n,n'$ etc.~will be determined later in \S\ref{resolution.tau}. 
Note the identity $(Q')' = [-q_2^3 \, q_0, -q_2^3 \, q_1, -q_2^4] = Q$. 
We have proved the following: 
\begin{Theorem} \sl 
If $[F]$ is a nonsingular point in $\He$, then so is $[\E_\HI(F)]$. 
The assignment 
\[ \He \setminus \B \lra \He \setminus \B, \quad 
[F] \lra [\E_\HI(F)] \] 
is an involutive automorphism. In particular $\He$ is isomorphic to 
its own dual variety. 
\end{Theorem} 

\section{The Jacobian ideal} \label{section.J}
Let $J = \I(\E_\HI(\F))$ denote the Jacobian ideal of $\HI$. 

\subsection{} Let 
\begin{equation} 
0 \la R/J \la R \la R(-17) \otimes S_5 \la E_1 \la E_2 \la \dots 
\label{res1.J} \end{equation} 
denote the equivariant minimal resolution of $J$, i.e., 
$E_i$ is the module of $i$-th syzygies. 
Apply $\text{Hom}_R(-,R)$ to~(\ref{res1.J}) and consider the complex 
\[ 0 \ra R \stackrel{\epsilon_0}{\lra} R(17) \otimes S_5 
\stackrel{\epsilon_1}{\lra} E_1^\vee \ra \dots \] 
Write $E_1^\vee$ as a direct sum 
\[ \bigoplus\limits_{r \ge 1} \, R(17+r) \otimes M_r, \] 
where each $M_r$ is a finite direct sum of irreducible $SL_2$-representations. 

By construction $\epsilon_0(1)= \E_\HI(\F)$. 
Let $S_p \subseteq M_r$ denote a direct summand, and consider the 
composite 
\[ \theta: R(17) \otimes S_5 \lra  R(17+r) \otimes M_r \lra 
R(17+r) \otimes S_p. \] 
It can be seen as a map $S_5 \lra S_p$ whose coefficients are degree $r$ forms 
in the coefficients of $\F$. Hence, $\theta$ corresponds to a covariant $\Theta$ 
(determined up to a constant) of degree $r$ and order (say) $q$, defining 
\[ S_5 \lra S_p, \quad G \lra (G,\Theta)_{\frac{1}{2} (5-p+q)}. \] 
Altogether, the identity $\theta \circ \epsilon_0=0$ translates into 
\[ (\E_\HI(\F),\Theta)_{\frac{1}{2} (5-p+q)}=0. \] 

\subsection{First syzygies of $J$} \label{firstsyzygies.J} 
We will enumerate some of the first syzygies of $J$ by hand, and then 
show {\sl a posteriori} that they are a complete list. Since a syzygy in a certain 
degree produces non-minimal syzygies in higher degrees, at each stage we should 
ensure that only `new' syzygies are included. 

\begin{enumerate} 
\item[(i)] If $\II$ is any invariant of $d$-ics, then 
$(\E_\II(\F),\F)_{d-1} =0$ (see~Corollary~\ref{corollary.EPhi} below), 
hence $S_2$ is a summand in $M_1$. 
\item[(ii)]The space $\A_{5,5}$ is $2$-dimensional, and spanned by 
$\vartheta_{33} \, \vartheta_{22}$ and $\vartheta_{40} \, \F$. 
By construction 
$\widetilde I = (\E_\HI,\vartheta_{33} \, \vartheta_{22})_5$ is an 
invariant of degree $22$ (possibly zero). Since $\zeta_{22,0}=1$, we must have 
$\widetilde I = \alpha \, \vartheta_{40} \, \HI$ for some 
$\alpha \in \rationals$. Define 
\begin{equation} \cU = \vartheta_{33} \, \vartheta_{22} - \alpha \, 
\vartheta_{40} \, \F,  \label{defn.U} \end{equation}  
so that $(\E_\HI,\cU)_5=0$. 

\noindent Claim: This syzygy cannot have arisen from the submodule 
$S_2 \subseteq M_1$. 

\demo Otherwise it would correspond to a nonzero morphism $S_2 \otimes R_4 \lra S_0$. 
However $R_4 \simeq S_4(S_5)$ contains 
no copies of $S_2$ (or equivalently, $\zeta_{4,2} =0$ for quintics), hence this is impossible. 
\item[(iii)]By an analogous reasoning, 
if $\Phi$ is any covariant of degree-order $(9,5)$, then 
\[ (\E_\HI,\Phi)_5 = \text{some degree $8$ invariant} \times \HI. \] 
Since $\A_{8,0}$ has $\{\vartheta_{40}^2, \vartheta_{80} \}$ as a basis, 
this would produce a syzygy of the form 
\begin{equation} (\E_\HI,\Phi - 
\beta \, \vartheta_{40}^2 \, \F - \gamma \, \vartheta_{80} \, \F)_5 =0 \quad 
\text{for some $\beta,\gamma \in \rationals$.} 
\label{syz.bc} \end{equation} 
However, we need to weed out those syzygies which come from earlier 
degrees. Broadly speaking, we have three syzygies in degree $9$ which arise 
in this way, amongst which two come from earlier degrees and one will be new. 
The space $\A_{9,5}$ is $5$-dimensional with a basis 
(see page~\pageref{basis.95})
\begin{equation} \begin{array}{rrrrr} 
\vartheta_{51} \, \vartheta_{22}^2, & 
\vartheta_{51} \, \vartheta_{44},  & 
\vartheta_{40} \, \vartheta_{33} \, \vartheta_{22}, & 
\vartheta_{40}^2 \, \F, & \vartheta_{80} \, \F. 
\end{array} \label{A95basis} \end{equation}

The one-dimensional space $\A_{8,2}$ is spanned by $\vartheta_{82}$. 
From part (i) we get the obvious identity $((\E_\HI,\F)_4,\vartheta_{82})_2 =0$, which 
can be rewritten as $(\E_\HI,(\F,\vartheta_{82})_1)_5 = 0$. This is best seen 
symbolically. Writing $\E = e_\ux^5,\F = f_\ux^5, \vartheta_{82} = t_\ux^2$, 
both compound transvectants evaluate to $(e \, f)^4 \, (e \, t) (f \, t)$. 
Now $(\F,\vartheta_{82})_1$ is the following linear combination of 
the basis in~(\ref{A95basis}): 
\[ \qquad 
-\frac{7}{10} \, \vartheta_{51} \, \vartheta_{22}^2 
-\frac{1}{4} \, \vartheta_{51} \, \vartheta_{44} 
+\frac{5}{12} \, \vartheta_{40} \, \vartheta_{33} \, \vartheta_{22} 
-\frac{1}{20} \, \vartheta_{40}^2 \, \F - \frac{1}{4} \, \vartheta_{80} \, \F. \] 
From (ii) we have the obvious syzygy $(\E_\HI,\vartheta_{40} \,  \cU)_5 = 0$. 
Let us define $\beta, \gamma \in \rationals$ such that the covariant 
\begin{equation} 
\cV = \vartheta_{51} \, \vartheta_{22}^2 - \beta \, \vartheta_{40}^2 \, \F - 
\gamma \, \vartheta_{80}\, \F \label{defn.V} \end{equation}
satisfies $(\E_\HI,\cV)_5 =0$. It is immediate that $\cV$ cannot be a linear 
combination of $(\F, \vartheta_{82})_1$ and $\vartheta_{40} \,  \cU$, hence we have a 
new syzygy. 
\end{enumerate} 
So far we have found three independent first syzygies of $J$ corresponding to 
$S_2 \subseteq M_1, S_0 \subseteq M_5, S_0 \subseteq M_9$. The rational 
numbers $\alpha,\beta,\gamma$ are uniquely determined by the identities 
$(\E_\HI,\cU)_5 = (\E_\HI, \cV)_5 =0$, but we do not yet know their values. 

\subsection{} We will now construct the morphism whose Hilbert-Burch 
complex is expected 
to give a resolution of $J$. Let us change our approach slightly, 
and let $\tau = (\alpha,\beta,\gamma)$ denote an {\sl arbitrary} triplet in $\rationals^3$. 
For $F \in S_5$, define 
\[ \begin{aligned} 
\sigma_{\tau}(F): \; & S_2 \oplus S_0 \oplus S_0 \lra S_5,  \\ 
& (A,c_1,c_2) \lra (A,F)_1 + c_1 \, \cU + c_2 \, \cV, 
\end{aligned} \] 
where $\cU, \cV$ are defined via formulae~(\ref{defn.U}),(\ref{defn.V}). 
Let $\Gamma_\tau$ denote the Wronskian of $\sigma_{\tau}(\F)$, which is a covariant 
of degree-order $(17,5)$. Let $\gb_\tau \subseteq R$ denote the ideal generated by the coefficients of $\Gamma_\tau$, 
and $V_\tau = V(\gb_\tau) \subseteq \P^5$ the corresponding subvariety.  One knows {\sl a priori} that each 
of the components of $V_\tau$ is of codimension $\le 2$. 
We claim that $(\Gamma_\tau,\F)_5 = \constant \, \HI$ for all $\tau$. Indeed, the left hand side 
is a degree $18$ invariant, hence a numerical multiple of $\HI$. It remains to 
check that it does not vanish identically, which is easily verified by specializing to 
$x_1^5 + x_2^5 + (x_1+2 \, x_2)^5$. It follows that $(\HI) \subseteq \gb_\tau$, 
hence $V_\tau \subseteq \He$. Since the latter contains no proper hypersurfaces, $\gb_\tau$ must be of pure height 
two. Hence the Eagon-Northcott complex (or what is the same, the Hilbert-Burch complex) of $\sigma_\tau$ is a 
minimal resolution of $\gb_\tau$. 

By a direct calculation, 
\begin{equation} \Gamma_\tau(\cF_Q) = - \frac{2^6 . 3^2 . 151 . 293}{5^{15}} \,  
q_2^3 \, (q_0 \, q_2 + 3 \, q_1^2) \, (5 \, q_0 \, q_2 - q_1^2) \, K_\tau \, \cF_{Q'},  
\label{Gamma.exp} \end{equation} 
where $K_\tau$ is the expression 
\begin{equation} \begin{aligned} 
{} & (75000 \, \gamma+28125) \, q_0^4 \, q_2^4+ \\ 
& (520000 \, \alpha+42000 \, \gamma-22500-960000 \, \beta) \, q_0^3 \, q_1^2 \, q_2^3+ \\ 
& (-1344000 \, \beta+292800 \, \gamma+872000 \, \alpha+6750) \, q_0^2 \, q_1^4 \, q_2^2+ \\ 
& (121200 \, \gamma-900-576000 \, \beta+408000 \, \alpha) \, q_0 \, q_1^6 \, q_2+ \\ 
& (12744 \, \gamma-69120 \, \beta+43200 \, \alpha+45) \, q_1^8. 
\end{aligned} \label{Ktau.exp} \end{equation} 
Since $\Gamma_\tau$ visibly vanishes on $\B$, we have $V_\tau \supseteq \B$. 
Now we would like to impose the condition that $V_\tau = \B$. This will happen iff 
$K_\tau$ is nonzero at every point of $\He \setminus \B$, i.e., iff 
\begin{equation} K_\tau = \delta \, q_2^r \, (q_0 \, q_2 + 3 \, q_1^2)^s \, 
(5 \, q_0 \, q_2 - q_1^2)^t  \label{Ktau.eqns} \end{equation} 
for some $\delta \in \rationals$, and nonnegative integers $r,s,t$ satisfying 
$r + 2s+2t =8$. It is easy to see that if we fix the choice of the triple 
$(r,s,t)$, then (\ref{Ktau.eqns}) is an inhomogeneous system of 
linear equations for the variables 
$\alpha,\beta,\gamma,\delta$. I solved this system in {\sc Maple}, and found that 
it admits a solution only in the following cases, the solution being unique in 
every case. 
\[ \begin{array}{ccc} 
(r,s,t) & & (\alpha,\beta,\gamma,\delta)\\ \hline 
(0,0,4) & & (0,0,0,1/45) \\ 
(0,2,2) & & (2/5,14/75,-2/5,-1/75) \\ 
(0,1,3) & & (1/6,2/45,-1/3,1/25).
\end{array} \] 
Thus we have the following theorem. 
\begin{Theorem} \sl 
\begin{enumerate} 
\item 
For any $\tau \in \rationals^3$, the ideal $\gb_\tau$ is perfect of height two with 
minimal resolution 
\[ \begin{aligned} 
0 \la R/\gb_\tau \la R & \la R(-17) \otimes S_5 \\ 
& \la R(-18) \otimes S_2 \oplus R(-22) \oplus R(-26) \la 0. 
\end{aligned} \] 
\item 
We have an inclusion of varieties $\B \subseteq V_\tau$, which is an 
equality iff $\tau$ is one of the following triples: 
\begin{equation}  
(0,0,0), \quad (2/5,14/75,-2/5) \quad (1/6,2/45,-1/3). 
\label{special.triples} \end{equation} 
\end{enumerate} 
\end{Theorem} 

\subsection{} \label{section.jacobian}
Now let $\tau = (\alpha,\beta,\gamma)$ be the {\sl specific} triple for which 
$(\E_\HI,\cU)_5 = (\E_\HI,\cV)_5 = 0$. We have shown that the complex 
\[ R \stackrel{\epsilon_0}{\lra} R(17) \otimes S_5 
\stackrel{\epsilon_1}{\lra} R(18) \otimes S_2 \oplus R(22) \oplus R(26) \] 
is exact in the middle. (Indeed, its middle cohomology is 
$\text{Ext}^1_R(R/\gb_\tau,R)$, which is zero since $\gb_\tau$ is perfect of 
height $2$.) Hence up to scalar, $\Gamma_\tau$ is the unique covariant of degree-order 
$(17,5)$ whose image by $\epsilon_1$ is zero. But $\E_\HI$ also has 
this property, hence $\Gamma_\tau = \text{nonzero constant} \times \E_\HI$. 
Since $V(J) = \B$, $\tau$ must be one of the special triples above, hence the following result: 
\begin{Proposition}  \sl 
For one of the three triples from~(\ref{special.triples}), we have the equality 
$\gb_\tau = J$ (the Jacobian ideal of $\HI$). In particular $J$ is also 
perfect. 
\end{Proposition} 
The value of $\tau$ will be found in~\S\ref{resolution.tau}. 
\subsection{The Cayley method} 
Initially I attempted to prove the perfection of $J$ by using the 
Cayley method of calculating resultants (see~\cite[Ch.~2]{GKZ}). This attempt failed, 
but the outcome was yet another perfect ideal supported on 
$\B$. Since the details are similar to~\cite[\S 5]{ACA}, we will be 
brief. 

Since $\HI$ is the resultant of $\F$ and $\vartheta_{33}$, it can be represented 
as the determinant of a complex. For a fixed $F \in S_5$, consider the Koszul complex 
\[ 0 \ra \O_{\P^1}(-8) \ra \O_{\P^1}(-3) \oplus \O_{\P^1}(-5) 
\stackrel{u}{\ra} \O_{\P^1} \ra 0, \] 
where $u$ is defined on the fibres as the map $(A,B) \ra A \, \vartheta_{33}(F) + B \, F$. 
Form the tensor product with $\O_{\P^1}(5)$, and consider the resulting hypercohomology spectral sequence. 
This produces a morphism 
\[ g_F: S_2 \oplus S_0 \oplus S_1 \lra S_5, \] 
such that $\det(g_\F) = \HI$. The component maps 
$S_2 \ra S_5, S_0 \ra S_5$ are easily described, they are 
$A \ra A \, \vartheta_{33}(F),    B \ra B \, F$. 
The third map $\mu_F: S_1 \ra S_5$ (which is a $d_2$-differential in the spectral sequence) 
is given via the {\sl Morley form},  described as follows: symbolically write $F = f_\ux^5, \vartheta_{33} = c_\ux^3$, 
and define 
\[ \M = (f \, c) \, [ \, f_\ux \, f_\uy^3 \, c_\uy^2 + 
c_\ux \, f_\uy^4 \, c_\uy \, ]. \] 
This defines a {\sl bivariate} covariant of $\F$,  of orders $1$ and $5$ respectively in $\ux,\uy$. 
If $A \in S_1$, then $\mu_F(A) = (A,\M)_1$. 
(The transvectant is with respect to $\ux$-variables, so the result is an order $5$ form in $\uy$.) 
We may instead decompose $\M$ into its Gordan series, and write 
\[ \mu_F(A) = (A,(F,\vartheta_{33})_1)_1 + 
\frac{1}{6} \, A \, (F,\vartheta_{33})_2. \] 
Now consider the truncated morphism $h_F: S_2 \oplus S_1 \lra S_5$, and let $\Lambda$ denote the 
Wronskian of $h_\F$. By construction $\Lambda$ is also a covariant of degree-order 
$(17,5)$, moreover $(\Lambda,\F)_5 = \constant \, \HI$. (Compare the argument in the previous section.) 
By a direct calculation, 
\[ \Lambda(\cF_Q) = 
- \frac{2^{16} \, 3^9}{5^{14}} \, q_2^3 \, 
(q_0 \, q_2 + 3 \, q_1^2)^2 \, (5 \, q_0 \, q_2 - q_1^2)^5 \, x_1^5, \] 
hence $\Lambda$ differs from $\E_\HI$ (or rather from any of the $\Gamma_\tau$). However, 
$\Lambda$ vanishes exactly over $\B$, hence by the usual argument we 
get the following result: 
\begin{Proposition} \sl 
The ideal $\I(\Lambda)$ is perfect of height two, with equivariant minimal resolution 
\[ 0 \la R/{\I(\Lambda)} \la R \la R(-17) \otimes S_5 \la 
R(-20) \otimes S_2 \oplus R(-21) \otimes S_1 \la 0. \] 
\qed \end{Proposition} 

\section{Evectants} \label{section.evectants} 
We could resolve the ambiguity about the correct value of $\tau$ (and hence about the maps in the 
resolution of $J$), if we could only derive an expression for $\E_\HI$. 
It is certainly possible to compute the latter in {\sc Maple} by a brute-force 
differentiation, but I have avoided this route in the belief that the general formalism developed 
here will prove useful elsewhere. 

In this section the construction of evectants will be generalized as follows: 
given any covariant $\Phi$ of $d$-ics we will associate to it a sequence of covariants called the evectants of 
$\Phi$. We will then deduce formulae for the evectants of $(\Phi,\Psi)_r$ in terms of those of $\Phi$ and $\Psi$. 
Finally this machinery will be applied to formula~(\ref{defn.H}). 
We will heavily use the symbolic method, however the final result of the calculation can be understood 
(and used) without any reference to it. Additional variable-pairs $\uy,\uz$ etc will be used as necessary, and then 
$\Omega_{\uy \uz}$ etc denote the corresponding Omega operators. 
\subsection{Evectants of a covariant} Let 
$\F = f_\ux^d = \sum\limits_{i=0}^d \, \binom{d}{i} \, a_i \, 
x_1^{d-i} \, x_2^i$
denote a generic binary $d$-ic, and let 
\[ \E(\ux) = \sum\limits_{i=0}^d \, 
\frac{\partial}{\partial a_i} \, 
x_1^i \, (-x_2)^{d-i} \, \] 
denote the evectant operator. 
Let $\Phi = \varphi_\ux^n$ be a covariant of degree-order $(m,n)$ of $d$-ics. Define 
\begin{equation} \Gamma = \frac{1}{m} \, [ \, \E(\ux) \circ \varphi_\uy^n \, ],  
\label{eqn1.gamma} \end{equation} 
which is a bihomogeneous form of orders $d,n$ in $\ux,\uy$ respectively, so that 
\[ (\Gamma,\F)_d = \frac{1}{m} \, 
\sum\limits_{i=0}^d \, a_i \, \frac{\partial \, \Phi(\uy)}{\partial a_i} = \Phi(\uy). \] 
Expanding $\Gamma$ into its Gordan series (\S\ref{section.Gordanseries}), we may write 
\begin{equation} \Gamma = \sum\limits_{i=0}^{\min(d,n)} \, 
(\ux \, \uy)^i \, {\alpha_{(i)}}_{\ux}^{d-i} \, 
{\alpha_{(i)}}_{\uy}^{n-i}, 
\label{eqn2.gamma} \end{equation} 
where $\A_i = {\alpha_{(i)}}_{\ux}^{d+n-2i}$ are a series of 
covariants of $f_\ux^d$. 
Now apply $(-,f_\ux^d)_d$ to each term in (\ref{eqn2.gamma}). Since 
\[ ((\ux \, \uy)^i \, \alpha_{\ux}^{d-i} \, 
\alpha_{\uy}^{n-i},f_\ux^d)_d = 
(\alpha \, f)^{d-i} \, \alpha_{\uy}^{n-i} \, f_\uy^i = 
[  (\alpha_\ux^{d+n-2i}, f_\ux^d)_{d-i} \, ]_{\ux:=\uy}, \] 
we deduce the identity 
\begin{equation} \sum\limits_i \, (\A_i,\F)_{d-i} = \Phi. 
\label{identity.Ai.F} \end{equation} 
The covariants $\A_\bullet = \{\A_0,\A_1,\dots,\A_{\min(d,n)}\}$ will 
be called the evectants of $\Phi$. By construction $\A_i$ is of 
degree-order $(m-1,d+n-2i)$.  
If $\Phi$ is an invariant, then $\A_0$ (the only nonzero 
evectant) coincides with $\E_\Phi$ as defined in \S\ref{section.evectant}. 

\begin{Lemma} \sl With notation as above, 
\[ \A_i = \frac{(d+n-2i+1)!}{i! \, (d+n-i+1)! \,m} 
\, \left\{ \Omega_{\ux \, \uy}^i  \circ 
[ \E(\ux) \circ \Phi(\uy)] \right\}_{\uy := \ux \, .} \] 
\label{lemma.Ev1} \end{Lemma} 
\demo Apply $\Omega_{\ux \uy}^\ell$ to each term in~(\ref{eqn2.gamma}), 
and use Lemma~\ref{lemma.Ev2} below. 
The terms with $\ell > i$ vanish because 
$\Omega_{\ux \uy} \circ \alpha_{\ux}^{d-i} \, \alpha_\uy^{n-i} =0$. 
Those with $\ell < i$ vanish after we set $\uy:=\ux$, this leaves only the term $\ell = i$. \qed 

\subsection{The evectants of a transvectant} 
Let $\Phi = \varphi_\ux^n, \, \Psi = \psi_\ux^{n'}$ denote two 
covariants with degree-orders $(m,n),(m',n')$, and evectants 
$\A_\bullet, \B_\bullet$ respectively. Their $r$-th transvectant 
$\Theta = (\Phi,\Psi)_r$ is of degree-order $(m+m',n+n'-2r)$. 
We would like to deduce formulae for the evectants $\C_\bullet$ of 
$\Theta$ in terms of the data $\Phi,\Psi,\A_\bullet,\B_\bullet$. Let us write 
\[ \Theta(\uy) = \frac{(n-r)! \, (n'-r)!}{n! \, n'!} \, 
\left \{ \Omega_{\uy \uz}^r \circ [ \Phi(\uy) \, \Psi(\uz) \, ] 
\right \}_{\uz:= \uy \, }    \] 
(if we expand $\Omega_{\uy \uz}^r$ by the binomial theorem, then this reduces to the 
definition in \S\ref{section.trans}), and then 
\[ \C_s = \kappa \, \left \{ \Omega_{\ux \uy}^s \circ \underbrace{[ 
\, \E(\ux) \circ \left \{ \Omega_{\uy \uz}^r \circ [ \Phi(\uy) \, \Psi(\uz) \, ] 
\right \}_{\uz:= \uy } ]}_{\en{a}} \right\}_{\uy:=\ux \, ,} \] 
where 
\begin{equation} 
\kappa = \frac{(n-r)! \, (n'-r)! \, (d+n+n'-2r-2s+1)!} 
{n! \, n'! \, s! \, (d+n+n'-2r-s+1)! \, (m+m')} \, .
\label{eqn.c1} \end{equation} 
It is understood that $r \le \min(n,n')$ and $s \le \min(d,n+n'-2r)$. 

The operators $\E(\ux)$ and $\Omega_{\uy \uz}$ commute, since they involve 
disjoint sets of variables. Hence 
\[ \en{a} = [ \, \Omega_{\uy \uz}^r \circ \underbrace{ 
\left\{ \, \E(\ux) \circ [\, \Phi(\uy) \, \Psi(\uz) \, ] \, 
\right\} }_{\en{b}} \;  ]_{\uz:=\uy \, .}  \] 
By the product rule for differentiation, 
\[ \en{b} = \underbrace{[ \, \E(\ux) \circ \Phi(\uy) \, ] \, \Psi(\uz)}_{\en{b_1}} + 
\underbrace{\Phi(\uy) \, [ \, \E(\ux) \circ \Psi(\uz) \, ]}_{\en{b_2}}. \] 
\subsection{} 
Writing $\A_i = {\alpha_{(i)}}_{\ux}^{d+n-2i}$, 
\begin{equation}  \en{b_1} = m \, \sum\limits_i \, 
(\ux \, \uy)^i \, {\alpha_{(i)}}_{\ux}^{d-i} \, 
{\alpha_{(i)}}_{\uy}^{n-i} \, \psi_\uz^{n'}.  \label{b.term1} \end{equation} 
We have to apply $\Omega_{\uy \uz}^r$ to each term in~$\en{b_1}$, 
and then set $\uz:=\uy$. The recipe is best seen combinatorially (also see~\cite[\S3.2.5]{Glenn}). 
From each summand in~(\ref{b.term1}) we sequentially remove $r$ symbolic 
factors involving $\uy$, and pair them with similarly removed $r$ factors involving $\uz$. 
By pairing a factor of the type $\beta_\uy$ with one of the type $\gamma_\uz$, we get a new factor 
$(\beta \, \gamma)$. 

The $\uz$-factors are all necessarily equal to $\psi_\uz$, on the other hand we may suppose that 
$k$ of the $\uy$-factors are $(\ux \, \uy)$ and the rest $r-k$ are ${\alpha_{(i)}}_{\uy}$. It is convenient 
to see $(\ux \, \uy)$ as $h_\uy$ with $(h_1,h_2) = (-x_2,x_1)$. Then the pairings produce 
factors $(h \, \psi)^k = (-1)^k \, \psi_\ux^k$ and $(\alpha_{(i)} \, \psi)^{r-k}$ respectively. 
We think of the $r$ copies of $\Omega_{\uy \uz}$ operating one after the 
other, so that the temporal sequence of removing the factors 
needs to be taken into account. At any stage, we may remove 
an $(\ux \, \uy),\psi_\uz$ pair or an ${\alpha_{(i)}}_{\uy},\psi_\uz$ pair, hence 
there are $r!/(k! \, (r-k)!)$ ways of choosing this sequence. 
The $\psi_\uz$ factors which have been removed can be sequentially ordered in 
$\frac{n'!}{(n'-r)!}$ ways (regarding them as mutually distinguishable), with 
a similar argument for other factors. This gives the expression 
\[ \begin{aligned} 
{} & [ \, \Omega_{\uy \uz}^r \circ \en{b_1} \,]_{\uz:=\uy} = \\ 
& m \, 
\sum\limits_i \, \sum\limits_k    \, \lambda(i,k;n,n') \; 
\underbrace{(\ux \, \uy)^{i-k} \, (\alpha_{(i)} \, \psi)^{r-k} \, 
{\alpha_{(i)}}_{\ux}^{d-i} \, 
{\alpha_{(i)}}_{\uy}^{n-i-r+k} \, \psi_\ux^k \;  \psi_\uy^{n'-r} \,}_{\en{c}}, 
\end{aligned} \] 
where 
\begin{equation} 
\lambda(i,k;n,n') = (-1)^k \, \frac{r!}{k! \, (r-k)!} \, 
\frac{i!}{(i-k)!} \, \frac{(n-i)!}{(n-i-r+k)!} \, \frac{n'!}{(n'-r)!}. 
\label{eqn.c2} \end{equation} 
The inner sum is quantified over 
$\max(0,r-n+i) \le k \le \min(i,r)$, which is exactly the 
possible range of removals. Our numerical assumptions imply that 
the range is always nonempty. 

The reader who dislikes the combinatorial argument may verify the formula 
\[ \Omega_{\uy \uz} \circ \alpha_\uy^p \, \beta_\uz^q = 
p \, q \,  (\alpha \, \beta) \, \alpha_\uy^{p-1} \, \beta_\uz^{q-1} \] 
by a direct calculation, and then proceed by induction. 

\subsection{} 
The next task is to apply $\Omega_{\ux \uy}^s$ to $\en{c}$, and 
then set $\uy:=\ux$. We need a preliminary lemma which describes how the 
operator $\Omega_{\ux \uy}$ can be 
`cancelled' against a factor of $(\ux \, \uy)$. 
\begin{Lemma} \sl For integers $p,q,\ell,i \ge 0$, 
we have an equality 
\[[\, \Omega_{\ux \uy}^\ell \circ 
 (\ux \, \uy)^i \, a_\ux^p \, b_\uy^q \, ]_{\uy:=\ux} = 
\begin{cases} 
\mu(p,q;\ell,i) \; 
[\, \Omega_{\ux \uy}^{\ell-i} \circ a_\ux^p \, b_\uy^q \, ]_{\uy:=\ux}  
& \text{if $\ell \ge i$}, \\ 
0 & \text{otherwise,} \end{cases} \] 
where 
\begin{equation} \mu(p,q;\ell,i) = \frac{\ell!}{(\ell-i)!} 
\frac{(p+q-\ell+2i+1)!}{(p+q-\ell+i+1)!} \, . 
\label{eqn.c3} \end{equation} 
\label{lemma.Ev2} \end{Lemma} 
\demo Let $\G$ denote an arbitrary bihomogeneous form of orders $p,q$ in $\ux,\uy$ respectively. 
By straightforward differentiation, 
\[ \begin{aligned} 
{} & \Omega_{\ux \uy}  \circ (\ux \, \uy) \, \G = 2 \, \G + 
 (x_1 \, \frac{\partial \G}{\partial x_1} +  
x_2 \, \frac{\partial \G}{\partial x_2}) +  
(y_1 \, \frac{\partial \G}{\partial y_1} +  
y_2 \, \frac{\partial \G}{\partial y_2}) + \\ 
& (x_1 \, y_2 - x_2 \, y_1) \, 
(\frac{\partial^2 \G}{\partial x_1 \, \partial y_2} - 
\frac{\partial^2 \G}{\partial x_2 \, \partial y_1}) \\ 
= \;  & (p+q+2) \, \G + (\ux \, \uy) \, \Omega_{\ux \uy} \circ \G. 
\end{aligned} \] 
Now proceed by induction on $\ell,i$, and observe that terms 
involving $(\ux \, \uy)$ vanish once we set $\uy:=\ux$. \qed 

\medskip 

\noindent Hence 
\[ \en{d} = 
[ \, \Omega_{\ux \uy}^s \circ \en{c} \, ]_{\uy:=\ux} \] 
vanishes if $s < i-k$. Assume $s \ge i-k$, then 
\begin{equation} \begin{aligned} 
\en{d} =  & \, \mu(d-i+k,n+n'-2r-i+k \, ;s,i-k) \, \times \\
& (\alpha_{(i)} \, \psi)^{r-k} \, 
\underbrace{ [\, \Omega_{\ux \uy}^{s-i+k} \circ 
{\alpha_{(i)}}_{\ux}^{d-i} \, 
{\alpha_{(i)}}_{\uy}^{n-i-r+k} \, 
\psi_\ux^k \;  \psi_\uy^{n'-r} \, ]_{\uy:=\ux}}_{\en{e}}. 
\end{aligned} \label{eqn.de} \end{equation} 
Now $\en{e}$ can be evaluated using the following lemma. 
\begin{Lemma} \sl For integers $p_1,q_1,p_2,q_2,u \ge 0$, 
we have an equality 
\[ \begin{aligned} 
{} & [ \, \Omega_{\ux \uy}^u \circ 
a_\ux^{p_1} \, a_\uy^{q_1} \, b_\ux^{p_2} \, 
b_\uy^{q_2} \, ]_{\uy:=\ux} \\ 
= \; & \nu(p_1,q_1,p_2,q_2;u) \times 
(a \, b)^u \, a_\ux^{p_1+q_1-u} \; b_\ux^{\, p_2+q_2-u}, 
\end{aligned} \] 
where $\nu$ is given by the sum 
\[ \sum\limits_t \, 
(-1)^{u-t} \, \frac{u!}{t! \, (u-t)!} \, 
\frac{p_1!}{(p_1-t)!} \, \frac{q_1!}{(q_1-u+t)!} \, 
\frac{p_2!}{(p_2-u+t)!} \, \frac{q_2!}{(q_2-t)!}, \] 
quantified over 
\[ \max(0,u-\min(q_1,p_2)) \le t \le 
\min(p_1,q_2,u). \] 
(The sum is understood to be zero if this range is 
empty.) 
\end{Lemma} 
\demo This is essentially the same combinatorial argument as before. 
Note however that since $(a \, a) =0$, we cannot 
pair $a_\ux$ with $a_\uy$, and similarly for $b$. 
Assume that we have removed respectively $t,u-t,u-t,t$ copies of 
$a_\ux,a_\uy,b_\ux,b_\uy$. Then pairings of $a_\ux,b_\uy$ produce 
$(a \, b)^t$, and those of $b_\ux,a_\uy$ produce 
$(b \, a)^{u-t} = (-1)^{u-t} \, (a \, b)^{u-t}$. 

The range of $t$ is exactly such that the removals are possible, 
e.g., $u-t$ cannot exceed $q_1$ or $p_2$ etc. \qed 

\medskip 

\noindent It follows that 
\[ \begin{aligned} 
(\alpha_{(i)} \, \psi)^{r-k} \, \en{e} =  
& \, \nu(d-i,n-i-r+k,k,n'-r;s-i+k) \, \times \\ 
& \underbrace{(\alpha_{(i)} \, \psi)^{r-i+s} \, 
{\alpha_{(i)}}_\ux^{d+n-r-s-i} \, \psi_\ux^{n'-r-s+i}}_{\en{f}}, 
\end{aligned} \] 
and of course $\en{f} = (\A_i,\Psi)_{r-i+s}$. The calculation for $\en{b_2}$ is 
essentially the same, hence we are done. 

\begin{Theorem} \sl With notation as above, 
\begin{equation} 
\C_s = \sum\limits_{i=0}^{\min(d,n)} \,  \xi_i \, (\A_i,\Psi)_{r-i+s} + 
\sum\limits_{i=0}^{\min(d,n')} \,  \eta_i \, (\B_i,\Phi)_{r-i+s}, 
\label{formula.Cs} \end{equation} 
where 
\[ \begin{aligned} 
\xi_i = \kappa \, m \, \sum\limits_k \; \{ \, 
& \lambda(i,k;n,n') \; \mu(d-i+k,n+n'-2r-i+k \, ;s,i-k) \, \times \\ 
& \nu(d-i,n-i-r+k,k,n'-r \, ;s-i+k) \, \}, \end{aligned} \] 
and 
\[ \begin{aligned} 
\eta_i = (-1)^r \, \kappa \, m' \, \sum\limits_k \; \{ \, 
& \lambda(i,k;n',n) \; \mu(d-i+k,n+n'-2r-i+k \, ;s,i-k) \, \times \\ 
& \nu(d-i,n'-i-r+k,k,n-r;s-i+k) \, \}. \end{aligned} \] 
The sums are respectively quantified over 
\[ \begin{aligned} 
{} & \max(0,r-n+i,i-s) & \le k \le \min(i,r), \\ 
& \max(0,r-n'+i,i-s) & \le k \le \min(i,r). 
\end{aligned} \] 
\end{Theorem} 
\subsection{} \label{resolution.tau}
Note the following classical proposition: 
\begin{Proposition} \sl A degree $m$ covariant $\Phi$ is a ${\mathbf Q}$-linear combination of 
compound transvectants 
\[ (\dots ((\F,\F)_{r_1},\F)_{r_2},\dots,\F)_{r_{m-1}}. \] 
\end{Proposition} 
\demo This is usually proved using the symbolic method, but 
it is easy to give an alternate proof. Write $\Phi = 
\sum\limits_i \, (\A_i,\F)_{d-i}$, and use the inductive hypothesis to write 
each $\A_i$ in terms of compound transvectants. \qed 

\smallskip 

The only nonzero evectant of $\Phi = \F$ is $\A_d = 1$. Starting from this, 
in principle we can calculate the evectants of any covariant. 
I have programmed formula~(\ref{formula.Cs}) in {\sc Maple}, so that 
the calculations can be made seamlessly. 

\begin{Example} \rm 
Let $d=5, \Phi = (\F,\F)_2, \Psi = (\F,\F)_4$ and 
$\Theta = (\Phi,\Psi)_1$. 
Now $\Phi,\Psi$ have only one nonzero evectant each, namely 
$\A_3 = \F, \B_1 = \F$. Hence $\Theta$ has evectants
\[ \begin{array}{ll} 
\C_0 = \frac{1}{4} \, \F \, \Phi, & 
\C_1 = \frac{2}{11} \, (\F,\Phi)_1, \\ 
\C_2 = - \frac{1}{4} \, \F \, \Psi 
- \frac{5}{18} \, (\F,\Phi)_2, & 
\C_3 = \frac{2}{7} \, (\F,\Psi)_1 
- \frac{10}{21} (\F,\Phi)_3, \\ 
\C_4 = \frac{3}{20} \, (\F,\Psi)_2 
- \frac{17}{56} \, (\F,\Phi)_4, &  
\C_5 = -\frac{2}{21} \, (\F,\Phi)_5. 
\end{array} \] 
In fact $\C_5$ vanishes identically, since quintics have no covariant of degree-order $(3,1)$.
\end{Example} 

It is now routine to calculate the evectants of $\vartheta_{22}^7$ and $ \F \, \vartheta_{39}$, and 
hence finally $\E_\HI$. The result is 
\begin{equation} \E_\HI(\cF_Q) = 
- \frac{2^6}{3.  \, 5^{14}} \; 
q_2^3 \, (q_0 \, q_2 + 3 \, q_1^2)^2 \, 
\, (5 \, q_0 \, q_2 - q_1^2)^4 \, \cF_{Q'}. 
\label{EH.expression} \end{equation}

\begin{Corollary} \sl We have an equality 
$J = \gb_{(1/6,2/45,-1/3)}$. 
\end{Corollary} 
\demo Comparing (\ref{Gamma.exp}) and (\ref{Ktau.eqns}) with (\ref{EH.expression}), 
we can read off the values $(r,s,t)=(0,1,3)$. \qed 

\medskip 

It is unnecessary to make six iterations in order to calculate the evectants of $\vartheta_{22}^7$. 
Instead observe that 
\[ 14 \, \Gamma = \E(\ux) \circ \vartheta_{22}(\uy)^7 = 7 \, \vartheta_{22}(\uy)^6 \, 
\underbrace{\E(\ux) \circ \vartheta_{22}(\uy)}_{(\ast)}, \] 
where $(\ast) = (\ux \, \uy) f_\ux^4 \, f_\uy$. Now one can find the Gordan series of $\Gamma$ directly. 

\subsection{} Let $\Phi$ be a covariant of degree-order $(m,n)$. The covariance property of $\Phi$ implies 
that its coefficients satisfy certain differential equations; this forces some identities between its 
evectants $\A_\bullet$. In this section we will make them explicit. 
Let $U = (u_0,u_1,u_1 \cb x_1,x_2)^2$ denote an arbitrary quadratic form. 
\begin{Proposition} \sl 
We have an equality 
\begin{equation} 
(( \, [\E(\ux) \circ \Phi(\uy)],\F(\ux) \, )_{d-1},
U(\ux) \, )_2 = 
\frac{n}{d} \, \{ \, (\Phi(\ux),U(\ux))_1 \}_{\ux:=\uy}. 
\label{Phi.diffeq} \end{equation}
\end{Proposition} 
By construction $\E(\ux) \circ \Phi(\uy)$ has orders $d,n$ in $\ux,\uy$. Its $(d-1)$-th transvectant 
with $\F(\ux)$ has $\ux$-order $2$, and finally the second transvectant with $U(\ux)$ has 
no $\ux$-variables remaining. Thus both sides are order $n$ forms in $\uy$. 

\begin{Corollary} \sl 
If $\Phi$ is an invariant, then 
$(\E_\Phi,\F)_{d-1} = 0$. 
\label{corollary.EPhi} \end{Corollary} 
\demo 
The right hand side of (\ref{Phi.diffeq}) vanishes, 
hence the second transvectant of an arbitrary quadratic form with 
$(\E_\Phi,\F)_{d-1}$ is zero. This forces the latter to be zero. \qed 

\smallskip 

We sketch a proof of the proposition. Since both sides are linear in $U$, it suffices to check the identity for 
each of the basis elements $\{x_1^2, x_1 \, x_2, x_2^2 \}$. After unravelling the transvectants, 
we are reduced to the following differential equations known to be 
satisfied by any covariant (see~\cite[\S1.2.12]{Glenn}): 
\begin{equation} \begin{aligned} 
{} & \sum\limits_{i=0}^{d-1} \, (d-i) \, a_{i+1} \, 
\frac{\partial \Phi}{\partial a_i}  = 
x_1 \, \frac{\partial \Phi}{\partial x_2}, \qquad 
\sum\limits_{i=1}^d i \, a_{i-1} \, \frac{\partial \Phi }{\partial a_i} = 
x_2 \, \frac{\partial \Phi}{\partial x_1}, \\ 
& \sum\limits_{i=0}^d \, (d-2i)  \, a_i \, \frac{\partial \Phi }{\partial a_i} = 
x_1 \, \frac{\partial \Phi}{\partial x_1} - 
x_2 \, \frac{\partial \Phi}{\partial x_2}. 
\end{aligned} \label{cayley.equations} 
\end{equation} \qed 

Broadly speaking, these equations express the fact that 
$\Phi$ is annihilated by the three generators 
$\left(\begin{array}{rr} 0 & 1 \\ 0 & 0 \end{array}\right), 
\left(\begin{array}{rr} 0 & 0 \\ 1 & 0 \end{array}\right), 
\left(\begin{array}{rr} 1 & 0 \\ 0 & -1 \end{array}\right)$ 
of the Lie algebra ${\mathfrak {sl}}_2$. 
\subsection{} 
It turns out that one can remove the reference to $U$ from~(\ref{Phi.diffeq}) and rephrase 
it as a set of three identities in the $\A_\bullet$. They will be of the form 
\[ 
\sum\limits_i \, \omega_{i,q} \, (\A_i,\F)_{d-i-1+q} = 0, \qquad (q=0,1,2), \] 
where $\omega_{i,q}$ are certain rational numbers. 
The calculations are thematically similar to those 
we have just seen, so the derivation will only be sketched. 

Write $\E(\ux) \circ \Phi(\uy) = p_\ux^d \, q_\uy^n, \, 
\F = f_\ux^d, \, U = u_\ux^2, \, \Phi = \varphi_\ux^n$. Then the left and right hand 
sides of (\ref{Phi.diffeq}) respectively equal 
\[( \, \underbrace{(p \, f)^{d-1} \, p_\ux \, f_\ux \, q_\uy^n}_{(\star)},
u_\ux^2)_2, \quad 
( \, \underbrace{\frac{n}{d} \, (\ux \, \uy) \, \varphi_\ux \, 
\varphi_\uy^{n-1}}_{(\star \star)},u_\ux^2)_2. \] 
The second transvectant of $Z = (\star) - (\star \star)$ with 
an arbitrary $U$ is zero, so $Z$ itself must be zero. Now substitute 
the sum $m \, \sum\limits_i \, (\ux \, \uy)^i \, {\alpha_{(i)}}_{\ux}^{d-i} \, 
{\alpha_{(i)}}_{\uy}^{n-i}$ for $\E(\ux) \circ \Phi(\uy)$, and expand $Z$ into 
its Gordan series. It is of the form 
\[ Z = \T_0 + (\ux \, \uy) \, \T_1 + 
(\ux \, \uy)^2 \, \T_2, \] 
where $\T_q$ are of orders $2-q,n-q$ in $\ux,\uy$. 
We get the required identities by writing $\T_q|_{\uy:=\ux}=0$. 
Although {\sl a priori} $\T_1$ involves $\Phi$, we can rewrite the 
latter in terms of $\A_\bullet$ using (\ref{identity.Ai.F}). 
In conclusion we have the following theorem: 
\begin{Theorem} \sl With notation as above, 
\begin{equation} 
\sum\limits_{i=0}^{\min(d,n)} \, 
\omega_{i,q} \, (\A_i,\F)_{d-i-1+q} = 0, 
\qquad (q=0,1,2), \label{identities.A.} 
\end{equation} 
where 
\[ \begin{array}{ll} 
\omega_{i,0} = d-i, & 
\omega_{i,1} = \frac{(d-i)(2i-n)}{d \, (n+2)}+
\frac{m \, i-n}{m \, d}, \\ 
\omega_{i,2} = i \, (n-i) \, (d-i+n+1). 
\end{array} \]
\end{Theorem} \qed 

These identities are collectively equivalent to (\ref{cayley.equations}), 
but in contrast to the latter, each of them is individually invariant under a change of co{\"o}rdinates. 

If a sequence of covariants $\{\A_i\}$ of degree-orders $(m,d+n-2i)$ 
is to appear as the sequence of evectants of some $\Phi$, it is necessary that they satisfy the 
identities above. It would be of interest to find a set of necessary and sufficient conditions. 
 
\bigskip 

{\sc Acknowledgements.} {\small 
This work was funded by the Natural Sciences and Engineering Research Council of Canada.
The following electronic libraries have been useful 
for accessing classical references:

\begin{itemize}
\item The G\"ottinger DigitalisierungsZentrum ({\bf GDZ}), 
\item Project Gutenberg ({\bf PG}), 
\item The University of Michigan Historical Mathematics 
Collection ({\bf UM}).
\end{itemize}}

{} 

\vspace*{3mm} 

\small 

\noindent {\sc Jaydeep Chipalkatti} \\ 
433 Machray Hall \\
Department of Mathematics \\ University of Manitoba \\
Winnipeg, MB R3T 2N2 \\ Canada. \\ 
e-mail: {\tt chipalka@cc.umanitoba.ca}


\begin{thebibliography}{10}

\bibitem{AC}
A.~Abdesselam and J.~Chipalkatti. 
\newblock On the Wronskian combinants of binary forms. 
\newblock To appear in \emph{J.~Pure and Appl.~Algebra}. 

\bibitem{Aluffi-Faber}
P.~Aluffi and C.~Faber.
\newblock Linear orbits of $d$-tuples of points in 
${\mathbb P}^1$.
\newblock {\em J.~Reine Angew.~Math.}, band 445, pp.~205--220, 1993.

\bibitem{ACA}
C.~D'Andrea and J.~Chipalkatti (with an appendix by 
A.~Abdesselam). 
\newblock On the Jacobian ideal of the binary discriminant. 
\newblock Preprint, 2006 (math.AG/0601705). 

\bibitem{ACGH}
E.~Arbarello, M.~Cornalba, P.A.~Griffiths and J.~Harris. 
\newblock {\em Geometry of Algebraic Curves, Volume~I}. 
\newblock Grundlehren der mathematischen Wissenschaften, No. 267, 
Springer-Verlag, 1985. 

\bibitem{Dolgachev}
I.~Dolgachev.
\newblock {\em Lectures on Invariant Theory}.
\newblock London Mathematical Society Lecture 
Notes No.~296. Cambridge University Press, 2003.

\bibitem{Ei} 
D.~Eisenbud. 
\newblock {\em Commutative Algebra, with a View Toward Algebraic Geometry}. 
\newblock Graduate Texts in Mathematics, Springer-Verlag, \, New York, 1995. 

\bibitem{FH}
W.~Fulton and J.~Harris.
\newblock {\em  Representation Theory, A First Course}.
\newblock Graduate Texts in Mathematics. Springer--Verlag,\,New York, 1991.

\bibitem{GKZ}
I.~Gelfand,  M.~Kapranov, and A.~Zelevinsky. 
\newblock  {\em Discriminants, Resultants and  Multidimensional     Determinants},
\newblock Birkh{\"a}user, Boston, 1994. 

\bibitem{Glenn}
O.~Glenn.
\newblock {\em The Theory of Invariants}.
\newblock Ginn and Co.,\,Boston, 1915 
({\bf PG}). 

\bibitem{Gordan}
P.~Gordan.
\newblock Das simultane System von zwei quadratischen quatern{\"a}ren Formen. 
\newblock Math.~Ann., vol.~56, pp.~1--48, 1903 ({\bf GDZ}). 

\bibitem{GrYo}
J.~H. Grace and A.~Young.
\newblock {\em The Algebra of Invariants, 1903.}
\newblock Reprinted by Chelsea Publishing Co., New York, 1962 ({\bf UM}). 

\bibitem{Gurevich}
G.~B. Gurevich. 
\newblock {\em Foundations of the Theory of Algebraic Invariants}. 
\newblock Noordhoff, Groningen, 1964. 

\bibitem{Harris}
J.~Harris. 
\newblock {\em Algebraic Geometry, a First Course}.
\newblock Graduate Texts in Mathematics. Springer--Verlag,\,New York, 1992.

\bibitem{Hermite1}
C.~Hermite.
\newblock Fonctions homog{\`e}nes a deux ind{\'e}termin{\'e}es.
\newblock {\em Cambridge and Dublin Math.~Journal}, May 1854. 
Reprinted in {\em {\OE}uvres de Charles Hermite}
(ed.~Emile Picard), Gauthier-Villars, Paris, 1905 ({\bf UM}). 

\bibitem{Hilbert1}
D.~Hilbert. 
\newblock {\"U}ber die Singularit{\"a}ten der Diskriminantenfl{\"a}che. 
\newblock Math.~Ann., vol.~30, pp.~437--441, 1887 ({\bf GDZ}). 

\bibitem{Iarrobino-Kanev}
A. Iarrobino and V. Kanev. 
\emph{Power Sums, Gorenstein Algebras and Determinantal Loci}. 
\newblock Springer Lecture Notes in Mathematics No.~1721, 1999. 

\bibitem{Kung-Rota}
J.P.S.~Kung and G.-C.~Rota.
\newblock The invariant theory of binary forms. 
\newblock Bulletin of the AMS., vol.~10, No.~1, pp.~27--85, 1984. 

\bibitem{Lang}
S.~Lang. 
\newblock {\em Algebra}. 
\newblock Addison-Wesley, 1984. 

\bibitem{MacMahon}
P.~MacMahon.
\newblock Article on `Algebraic Forms' \rm 
in the {E}ncyclop{\ae}dia {B}rittanica (11th edition), 1911.

\bibitem{Meulien}
M.~Meulien. 
\newblock Sur les invariants des pinceaux de formes quintiques binaires. 
\newblock Ann.~Inst.~Fourier (Grenoble), vol.~54, No.~1, pp.~21--51, 2004. 

\bibitem{Olver}
P.~Olver.
\newblock {\em Classical Invariant Theory}.
\newblock London Mathematical Society Student Texts. Cambridge University
  Press, 1999.

\bibitem{Salmon1}
G.~Salmon.
\newblock {\em Higher Algebra}. 
\newblock Reprinted by Chelsea Publishing Co., New York, 1964.

\bibitem{Sturmfels}
B.~Sturmfels.
\newblock {\em Algorithms in Invariant Theory}.
\newblock Texts and Monographs in Symbolic Computation. Springer--Verlag, Wien  New York, 1993.

\bibitem{Weyman}
J.~Weyman. 
\newblock The equations of strata for binary forms. 
\newblock J.~Algebra, vol.~122, pp.~244--249, 1989. 
\end{thebibliography}
\end{document}